\documentclass{amsart}
\usepackage{amsmath,amsthm,bbm}
\usepackage{amsfonts,amssymb}

 \newtheorem{theorem}{Theorem}[section]
 \newtheorem{lemma}[theorem]{Lemma}
 \newtheorem{corollary}[theorem]{Corollary}
 \newtheorem{proposition}[theorem]{Proposition}
 \newtheorem{definition}[theorem]{Definition}
 
 \newtheorem{remark}[theorem]{Remark}
 
 \newcounter{figures}[section]

\def\bC{{\mathbb C}}

\def\NN{{\mathbb N}}

\def\RR{{\mathbb R}}
\def\SS{{\mathbb S}}
\def\ZZ{{\mathbb Z}}

  \def\cB{\mathcal{B}}
  \def\cC{\mathcal{C}}
  \def\cD{\mathcal{D}}

  \def\cI{\mathcal{I}}

  \def\cL{\mathcal{L}}
  \def\cM{\mathcal{M}}
  \def\cN{\mathcal{N}}
  \def\cO{\mathcal{O}}
  
  \def\cQ{\mathcal{Q}}
  \def\cR{\mathcal{R}}
  \def\cS{\mathcal{S}}

  \def\cX{\mathcal{X}}

\def\CB{{\mathcal B}}

\def\CV{{\mathcal V}}

\def\proj{\operatorname{proj}}
\def\Proj{\operatorname{Proj}}

\def\supp{\operatorname{supp}}

\def\eps{{\varepsilon}}
\def\a{{\alpha}}
\def\b{{\beta}}

\def\ha{{\hat a}}

\def\hc{{\hat c}}

\def\ph{{\varphi}}

\def\ONE{{\mathbbm 1}}
\def\tONE{{\tilde{\mathbbm 1}}}

\def\Fsrpq{{F_{pq}^{s\r}}}
\def\fsrpq{{f_{pq}^{s\r}}}

\def\Bsrpq{{B_{pq}^{s\r}}}

\def\bsrpq{{b_{pq}^{s\r}}}

\def\EE{E}
\def\cc{\tilde{c}}

\def\Q{{\mathcal Q}}

\def\tP{\tilde{P}}
\def\hA{{\hat A}}
\def\hB{{\hat B}}
\def\hC{{\hat C}}

\def\bargamma{{\tilde\gamma}}
\def\PsiZ{{F}}
\def\EE{{\cM}}

\def\r{{\rho}}

\def\WW{{W_{\a,\b}}}
\def\WWW{{W_{\mu,\a,\b}}}
\def\Lp{{p}}

\newcommand\F[4]{F^{#1#2}_{#3#4}}
\newcommand\B[4]{B^{#1#2}_{#3#4}}

\def\Fsrpq{{F_{pq}^{s\r}}}

\def\fsrpq{{f_{pq}^{s\r}}}

\def\Bsrpq{{B_{pq}^{s\r}}}

\def\bsrpq{{b_{pq}^{s\r}}}

\newcommand\norm[1]{\|#1\|}
\newcommand\nnorm[1]{\Big\|#1\Big\|}

\def\ip#1{{\langle #1 \rangle}}
\newcommand{\X}{\mathcal{X}}

\def\Lpp{{L^p(w_{\a,\b})}}

\def\Linfty{{\infty}}
\def\bC{{\mathbb C}}

\def\ph{{\varphi}}

\def\cc{{\rm c}}

\newcommand{\w}{w}

\def\cN{\mathcal{N}}
\def\NN{{\mathbb N}}

\def\LL2{{L^2(\w)}}
\def\tPnu{\tilde P_\nu^{(\a,\b)}}
\def\y{{\bar y}}

\begin{document}

\title[spaces of distributions induced by tensor product bases]
{Decomposition of spaces of distributions induced by tensor product bases}

\author{Kamen Ivanov, Pencho Petrushev and Yuan Xu}

\address{Institute of Mathematics and Informatics\\
Bulgarian Academy of Sciences\\ 1113 Sofia\\ Bulgaria}
\email{ivanov@math.sc.edu}

\address{Department of Mathematics\\University of South Carolina\\
Columbia, SC 29208\\
and Institute of Mathematics and Informatics, Bulgarian Academy of Sciences}
\email{pencho@math.sc.edu}

\address{Department of Mathematics\\ University of Oregon\\
Eugene, Oregon 97403-1222.}\email{yuan@math.uoregon.edu}

\date{February 9, 2009}



\thanks{The second author has been supported by NSF Grant DMS-0709046
and the thid author by NSF Grant DMS-0604056.}

\thanks{Address for manuscript correspondence:
Pencho Petrushev,
Dept. of Mathematics, University of South Carolina,
Columbia, SC 29208,
E-mail: pencho@math.sc.edu,
Phone: (803)777-6686.
}

\begin{abstract}
Rapidly decaying kernels and frames (needlets) in the context of tensor product
Jacobi polynomials are developed based on several constructions of multivariate
$C^\infty$ cutoff functions.
These tools are further employed to the development of
the theory of weighted Triebel-Lizorkin and Besov spaces on $[-1, 1]^d$.
It is also shown how kernels induced by cross product bases can be constructed and utilized for
the development of weighted spaces of distributions on products of multidimensional ball, cube,
sphere or other domains.
\end{abstract}

\maketitle

\pagestyle{myheadings}
\thispagestyle{plain}

\section{Introduction}\label{introduction}
\setcounter{equation}{0}

The purpose of this article is to introduce and study Triebel-Lizorkin and Besov spaces
on the d-dimensional cube $Q^d=[-1,1]^d$ with Jacobi weights
and discuss the respective spaces on the product domains $B^{d_2}\times Q^{d_1}$, $B^{d_1}\times B^{d_2}$
with $B^d$ being the unit ball in $\RR^d$
as well as sets of the form
$Q^{d_1}\times \SS^{d_2}$, $\SS^{d_1}\times B^{d_2}$, $Q^{d_1}\times T^{d_2}$,
$Q^{d_1}\times \RR^{d_2}$, $B^{d_1}\times \RR^{d_2}_+$,
$Q^{d_1}\times B^{d_2}\times\RR^{d_3}_+$, etc.
Here $\SS^{d}$ is the unit sphere in $\RR^{d+1}$,
$T^{d}$ is the simplex in $\RR^d$, and
$\RR^{d}_+:=[0, \infty)^d$.
In short, we are interested in developing the theory of distribution spaces
on some products of
$Q^{d_1}$, $B^{d_2}$, $\SS^{d_3}$, $T^{d_4}$, $\RR^{d_5}$, or $\RR^{d_6}_+$
with weights.
There are two important components of such undertaking:
(i) the spaces need to be properly defined and
(ii) building blocks need to be constructed and used for characterization of the spaces.
We maintain that for both tasks tensor product orthogonal bases should be used.

\subsection{The principle distinction between the spaces on \boldmath $[-1,1]^d$ and $B^d$}

It~seems to us natural to introduce weighted smoothness spaces on $[-1,1]^d$ or
$B^d$ with weights by means of orthogonal polynomials.
However, there is a surprising difference between the orthogonal polynomial expansions
on $[-1,1]^d$ and $B^d$ which we would like to describe next.

Let us first briefly review the definition of Triebel-Lizorkin and Besov spaces
on $B^d$, given in \cite{KPX2}.
Denote by $\CV_n^d$ the space of all polynomials of total degree $n$ which
are orthogonal to lower degree polynomials in $L^2(B^d, w_\mu)$
with weight
$
w_\mu(x):= (1-\|x\|_2^2)^{\mu-1/2}.
$
The orthogonal projector
$\Proj_n: L^2(B^d, w_\mu) \mapsto\CV_n^d$
can be written in the form
\begin{equation}\label{Proj-Vn}
(\Proj_n f)(x) = \int_{B^d} f(y)P_n(w_\mu; x,y) w_\mu(y) dy.
\end{equation}
To introduce weighted Triebel-Lizorkin ($F$-spaces) and Besov spaces ($B$-spaces) on $B^d$
(see \cite{Pee}, \cite{T1} for the general idea), let
\begin{equation}\label{def-kernels-sp}
\Phi_0(x, y) := 1
\quad\mbox{and}\quad
\Phi_j(x, y) := \sum_{n=0}^\infty \ha
\Big(\frac{n}{2^{j-1}}\Big)P_n(w_\mu; x,y), \quad j\ge 1,
\end{equation}
where $\ha\in C^\infty[0, \infty)$ is a cutoff function
such that $\supp \ha \subset [\frac{1}{2}, 2]$ and
$|\ha|\ge c >0$ on $[3/5, 5/3]$.

The weighted $F$-space $F^{s, q}_p$ on $B^d$ with $s\in \RR$, $0<p<\infty$, $0<q\le\infty$,
is defined as the space of all distributions $f$ on $B^d$ such that
\begin{equation}\label{Tri-Liz-norm0}
\|f\|_{F^{s, q}_p}
:=\Big\|\Big(\sum_{j=0}^{\infty}\big(2^{sj}|\Phi_j\ast f(\cdot)|\big)^q\Big)^{1/q}\Big\|_{L^p(w_\mu)}
<\infty,
\end{equation}
where $\Phi_j*f(x):= \langle f, \overline{\Phi_j(x, \cdot)}\rangle$
(as in \eqref{convolution1}).
The corresponding scale of weighted Besov spaces $B^{s, q}_p$
is defined via the (quasi-)norms
\begin{equation}\label{Besov-norm0}
\|f\|_{B^{s, q}_p}
:= \Big(\sum_{j=0}^\infty \Big(2^{s j}\|\Phi_j*f(\cdot)\|_{L^p(w_\mu)}\Big)^q\Big)^{1/q}.
\end{equation}
We refer the reader to \cite{KPX2} for more detailed account of weighed
F- and B-spaces on the ball.


A ``natural" attempt to introduce Triebel-Lizorkin and Besov spaces on $[-1, 1]^d$ with
weight
\begin{equation}\label{Jacobi-weight}
w_{\a,\b}(x):=\prod_{i=1}^d (1-x_i)^{\a_i}(1+x_i)^{\b_i}
\end{equation}
would be to use directly the same idea as above.
Namely, for multi-indices $\a$, $\b$, $\nu$
the d-dimensional tensor product Jacobi polynomials are defined by
\begin{equation} \label{product-Jacobi-intr}
\tP_\nu^{(\a, \b)}(x):=\prod_{j=1}^d \tP_{\nu_j}^{(\a_j, \b_j)}(x_j).
\end{equation}
Set
$\tP_n^{(\a, \b)}(x, y):= \sum_{|\nu|=n}\tP_\nu^{(\a, \b)}(x)\tP_\nu^{(\a, \b)}(y)$
and with $\ha$ as in (\ref{def-kernels-sp}) define
\begin{equation}\label{def-kernels-Q}
\Phi_0(x, y) := 1
\quad\mbox{and}\quad
\Phi_j(x, y) := \sum_{n=0}^\infty \ha
\Big(\frac{n}{2^{j-1}}\Big)\tP_n^{(\a, \b)}(x, y), \quad j\ge 1,
\end{equation}
which can be viewed as an analogue of the kernels from \eqref{def-kernels-sp}.

The next step would be to define weighted Triebel-Lizorkin and Besov spaces on $[-1, 1]^d$
with weight $w_{\a,\b}(x)$ exactly as in (\ref{Tri-Liz-norm0}) and (\ref{Besov-norm0})
using the kernels $\Phi_j(x, y)$ from (\ref{def-kernels-Q}).
Such a definition, however, is completely unacceptable due to the {\bf poor localization}
of the kernels $\Phi_j(x, y)$ from (\ref{def-kernels-Q}).
As~is shown in \cite{IPX} in the particular case of Legendre or Chebyshev polynomials, kernels of
the form (\ref{def-kernels-Q}) have no localization
whatsoever for some points $x, y \in[-1, 1]^d$.
In contrast, the kernels $\Phi_j(x, y)$ from (\ref{def-kernels-sp}) decay rapidly
away from the main diagonal in $B^d\times B^d$.
Interestingly enough, the situation is quite the same on the interval \cite{KPX1},
sphere \cite{NPW2}, simplex \cite{IPX}, and more surprisingly in the context of
tensor product Hermite \cite{PX3} and Laguerre functions \cite{KPPX}.

\subsection{The remedy for the problem}\label{remedy}
It appears that the tensor product Jacobi polynomials are in a sense of a different nature
compared to orthogonal polynomials on the interval, ball or simplex
as well as spherical harmonics and tensor product Hermite and Laguerre functions.
Truly multivariate cutoff functions need to be employed.
Our primary goal in this paper is to identify a natural class of cutoff functions
which will enable us to develop a meaningful theory of Triebel-Lizorkin and Besov
spaces on $[-1,1]^d$ with weight $w_{\a,\b}(x)$ via tensor product Jacobi polynomials.

The key is to consider multivariate cutoff functions $\hA$ with dyadic dilations covering
the whole spectrum and such that the kernels
\begin{equation}\label{def-kernels-introd}
\Phi_j(x, y) := \sum_{\nu\in\NN_0^d}
\hA\Big(\frac{\nu}{2^{j-1}}\Big)\tP_\nu^{(\a, \b)}(x)\tP_\nu^{(\a, \b)}(y)
\end{equation}
decay rapidly away from the main diagonal in $[-1, 1]^d\times [-1, 1]^d$.
It turns out that it suffices to consider compactly supported $C^\infty$ cutoff functions
$\hA:[0, \infty)^d\mapsto \bC$ which obey the following

\smallskip

\noindent
{\bf First Boundary Condition.}
{\em For any $t\in [0, \infty)^d$ which belongs to the coordinate planes, i.e.
$t=(t_1, \dots, t_{k-1}, 0, t_{k+1}, \dots, t_d)$ for some $1\le k \le d$,
\begin{equation}\label{first-kind}
\frac{\partial^m}{\partial t_k^m}\hA(t)=0,
\quad \mbox{for $\; m=1, 2, \dots$.}
\end{equation}
}
Sometimes, instead of this condition it will be more convenient to use the following
slightly more restrictive but for certain purposes better and easier to deal with

\smallskip

\noindent
{\bf Second Boundary Condition.}
{\em
There exists a constant $c_*>0$ such that for any $\tau\in [0, \infty)^d$
of the form $\tau=(\tau_1,\dots,\tau_{k-1},0,\tau_{k+1},\dots,\tau_d)$, $1\le k\le d$,
$\hA(t)= {\rm constant}$ for $t\in [\tau, \tau+c_*e_k]$
with  $e_k$ being the $k$th coordinate vector.
}

\smallskip

The point is that either of these conditions combined with $\hA$ being $C^\infty$
and compactly supported yields the rapid decay of the kernels $\Phi_j(x, y)$
from \eqref{def-kernels-introd} (see Theorem~\ref{thm:product-Jacobi-0}).
Then these kernels can be deployed to the definition of
weighted Triebel-Lizorkin and Besov spaces on $[-1, 1]^d$ by means of norms
similar to the norms in \eqref{Tri-Liz-norm0}-\eqref{Besov-norm0}.

As will be seen the weights
\begin{equation}\label{def-Wab-Intr}
W_{\a,\b}(n; x)
:= \prod_{i=1}^d (1-x_i + n^{-2})^{\a_i+1/2} (1+x_i + n^{-2})^{\b_i+1/2}
\end{equation}
will appear naturally in most estimates and results related to spaces on $[-1, 1]^d$
with weight $w_{\a,\b}(x)$.
Moreover, the inhomogeneity created by $w_{\a,\b}(x)$ and the boundary of $[-1, 1]^d$
leads us to the introduction
via $W_{\a,\b}(\cdot; \cdot)$ of a fourth parameter $\rho$ in the definition of
weighted Triebel-Lizorkin and Besov spaces on $[-1, 1]^d$.
Thus we introduce $F$-spaces by the norms
$$
\|f\|_{\Fsrpq}:=\Big\|\Big(\sum_{j=0}^{\infty}
\Big[2^{sj}\WW(2^j;\cdot)^{-\r/d}|\Phi_j*f(\cdot)|\Big]^q\Big)^{1/q}\Big\|_{L^p(w_{\a,\b})}
$$
and $B$-spaces by the norms
$$
\|f\|_{\Bsrpq} :=
\Big(\sum_{j=0}^\infty \Big[2^{s j}
\|\WW(2^j; \cdot)^{-\r/d}\Phi_j*f(\cdot)\|_{L^p(w_{\a,\b})}\Big]^q\Big)^{1/q}
$$
(see \S\S\ref{Tri-Liz}-\ref{Besov}).
This allows to use for different purposes various scales of
weighted $F$- and $B$-spaces on $[-1, 1]^d$.
For instance, as will be seen the Besov spaces $B_{\tau\tau}^{ss}$
appear naturally in nonlinear approximation in $L^p(w_{\a,\b})$.

As a next step we use kernels of the form \eqref{def-kernels-introd}
for the construction of building blocks (needlets) $\{\varphi_\xi\}$, $\{\psi_\xi\}$.
These are multiscale dual frames which enable us to characterize
the $F$- and $B$-space norms by the size
of the needlet coefficients $\{\langle f, \varphi_\xi \rangle\}$
in appropriate sequence norms.
They can be viewed as an analogue of the $\varphi$-transform of Frazier and Jawerth
\cite{FJ1, FJ2, FJW}.

The theory of weighted Triebel-Lizorkin and Besov spaces on $[-1, 1]^d$
and needlet decompositions in dimensions $d>1$ run parallel to
their theory in dimension $d=1$, developed in \cite{KPX1}, and on the ball \cite{KPX2}.
Therefore, to spare the reader the repetition of well established arguments
we shall only exhibit the essential differences and refer for the rest
to \cite{KPX1, KPX2}.
We shall place the emphasis on the development of multivariable cutoff functions
and the associated tensor product Jacobi kernels which defer substantially
from the ones in the univariate case and are the main reason for writing this paper.
We shall also discuss the main points of the development of Triebel-Lizorkin and Besov spaces
and needlets on products of
$[-1,1]^{d_1}$, $B^{d_2}$, $\SS^{d_3}$, $\RR^{d_4}$, or $\RR^{d_5}_+$
with weights as mentioned above.

This paper is part of a broader undertaking for development of spaces of distributions
in nonstandard settings such as on the sphere \cite{NPW2}, ball \cite{KPX2}
as well as in the frameworks of Hermite \cite{PX3} and Laguerre \cite {KPPX} expansions.
It is also closely related to the development of sub-exponentially localized Jacobi
and other kernels and needlets in \cite{IPX}.

\subsection{Outline of the paper}

A substantial part of the paper is devoted to the development of multivariate cutoff
functions and related tensor product Jacobi kernels.
In \S\ref{one-d-kernels}
we review some basic results from \cite{IPX, KPX1} and prove new results about
admissible univariate cutoff functions and the localization of the respective kernels
induced by univariate Jacobi polynomials.
In \S\ref{multi-kernels}
we present several constructions of multivariate admissible cutoff
functions.
We also construct cutoff functions of ``small" derivatives which
enables us to develop tensor product kernels with sub-exponential localization.
In \S\ref{d-Jacobi-kernels} the localization results of the corresponding
tensor product Jacobi polynomial are established.
In \S\ref{background}
we give some auxiliary results concerning a maximal operator and distributions
on $[-1, 1]^d$. We also establish some $L^p$-multipliers for tensor product
Jacobi polynomial expansions.
In \S\ref{Needlets}
we utilize kernels associated to cutoff functions of type (b) and (c) to
the construction of frame elements (needlets).
In \S\S\ref{Tri-Liz}-\ref{Besov}
we further use these kernels to define ``correctly" the weighted Triebel-Lizorlin
and Besov spaces on $[-1, 1]^d$ with weight $w_{\a,\b}(x)$.
We also establish needlet decomposition of the $F$- and $B$-spaces.
Section \ref{Nonlin-app} is devoted to nonlinear approximation from Jacobi needlets.
In \S\ref{cross-products} we briefly consider weighted spaces of distributions on $B^{d_1}\times [-1, 1]^d$.
In \S\ref{discussion} we discuss various aspects of distribution spaces on product
domains and tensor product bases.
Section \ref{appendix} is an appendix, where we place the lengthy proof of a lemma from \S\ref{cross-products}.

\smallskip

\noindent
{\bf Some useful notation.}
Throughout we shall denote
$$
  \|f\|_\Lp:=\Big(\int_{[-1,1]^d} |f(x)|^pw_{\a,\b}(x)dx\Big)^{1/p},
\quad 0<p<\infty,
$$
and
$
\|f\|_\Linfty:=\sup_{x\in[-1,1]^d}|f(x)|.
$
For $x\in\RR^d$ we shall use the norms
$\|x\|=\|x\|_\infty:=\max_i |x_i|$, $\|x\|_2:= (\sum_i |x_i|^2)^{1/2}$, and
$|x|=\|x\|_1:=\sum_i |x_i|$.
$\Pi_n^d$ will denote the set of all algebraic polynomials of total degree $\le n$ in $d$ variables.
Positive constants will be denoted by $c$, $c_1$, $c'$, $\dots$ and they may vary
at every occurrence, $a\sim b$ will stand for $c_1a\le b\le c_2a$.

\section{Localized Jacobi kernels induced by univariate cutoff functions}\label{one-d-kernels}
\setcounter{equation}{0}

Here we introduce the notion of admissible univariate cutoff functions
and review the localization properties of the associated kernels induced
by Jacobi polynomials established in \cite{IPX, KPX1, PX1}.
We also obtain some new localization results.

\subsection{Admissible univariate cutoff functions}
\label{principle}


\begin{definition}\label{cutoff-d1}
A function $\ha\in C^\infty [0, \infty)$ 
is said to be admissible if
 $\supp \ha \subset [0, 2]$ and $\ha^{(m)}(0)=0$ for $m\ge 1$.
Furthermore, $\ha$ is said to be admissible of
type $(a)$, $(b)$ or $(c)$ if $\ha$ is admissible and in addition obeys the respective condition:

$(a)$ $\ha(t) = 1$, $t\in [0, 1]$,

$(b)$ $\supp \ha \subset [1/2, 2]$ or

$(c)$ $\supp \ha \subset [1/2, 2]$ and $\sum_{j=0}^\infty |\ha(2^{-j}t)|^2=1$ for $t\in [1,\infty)$.
\end{definition}

We next introduce sets of $C^\infty$ functions with ``small'' derivatives.
As a tool for measuring the derivatives' growth we use functions $\cL$ satisfying the conditions:
\begin{equation}\label{def-M}
\begin{split}
\cL: [0, \infty)\mapsto[1, \infty)~\mbox{is monotone increasing,}~\cL(0)=1~\mbox{and}\\
M=M(\cL):= 1+\int_0^\infty \frac{dt}{(t+1)\cL(t)} <\infty.
\end{split}
\end{equation}
Typical examples of functions $\cL$ satisfying \eqref{def-M}
are
$\cL_{0,\eps}(t):=(1+t)^\eps$, $\varepsilon>0$, and
\begin{multline}\label{log1}
\cL_{\ell,\varepsilon}(t):=\ln(e+ t)\cdots
\underbrace{\ln\cdots\ln}_{\ell-1} \big(\underbrace{\exp\cdots\exp}_{\ell-1} 1+t\big)\\
\hspace{1.4in}\times
\big[\underbrace{\ln\cdots\ln}_{\ell}\big(\underbrace{\exp\cdots\exp}_{\ell}1+t\big)\big]^{1+\eps},
\end{multline}
where $\ell\in\NN$ and $0<\varepsilon\le1$.
Evidently,
$M(\cL_{\ell,\varepsilon}) \le c(\ell)\eps^{-1}$.

We shall use the standard notation $D_j^k:=\frac{\partial^k}{\partial x_j^k}$. 


\begin{definition}\label{Def-S}
Let $\cL$ satisfy \eqref{def-M}.
Given constants $\gamma,\bargamma>0$ and $d\ge 1$, we define $\cS(d, \cL;\gamma,\bargamma)$
to be the set of all functions
$\hA\in C^\infty[0,\infty)^d$, such that
$\|\hA\|_\infty \le \gamma$ and
\begin{equation}\label{der-L}
\frac{1}{k!}\|D_j^k \hA\|_\infty
\le \gamma\big(\bargamma\cL(k-1)\big)^k,
\quad\quad \forall k\in\NN,~1\le j\le d.
\end{equation}
\end{definition}

\smallskip

The next statement asserts the existence of admissible univariate cutoff functions with
``small'' derivatives.


\begin{theorem}\label{thm:cutoff-d=1}
Let $\cL$ and  $M$ be given by \eqref{def-M}.
Then the sets
$\cS(1,\cL; 1, 2M)$, $\cS(1,\cL; 2, 4M)$ and $\cS(1,\cL; 8, 8M)$
contain admissible cutoff functions $\ha$ of types
$(a)$, $(b)$, and $(c)$, respectively, $($see Definition~\ref{cutoff-d1}$)$ with values in $[0,1]$.
\end{theorem}
%

\noindent
{\bf Proof.}
We shall proceed quite similarly as in the proof of Theorem~3.1 in \cite{IPX}.
We~let
$\chi_\delta:= \frac{1}{2\delta}\ONE_{[-\delta, \delta]}$
and select
$\delta_j:=\frac{1}{(j+1)\cL(j)}$ for $j\ge 0$.
Apparently
$$
\sum_{j=0}^\infty \delta_j
\le 1+\int_0^\infty\frac{dt}{(t+1)\cL(t)}
= M.
$$
We define
$$
\varphi_m:=\chi_{\delta_0}* \dots * \chi_{\delta_m}
\quad\mbox{and}\quad
\varphi(t):=\lim_{m\to\infty} \varphi_m(t).
$$
Just as in \cite[Theorem 1.3.5]{H} we have
$\varphi\in C^\infty$, $\varphi\ge 0$, $\supp \varphi\subset [-M, M]$ and
$$
\|\varphi^{(k-1)}\|_\infty \le \frac{1}{\prod_{j=0}^{k-1} \delta_j}
\le k!\cL(k-1)^k
\quad\mbox{for}\quad k\ge 1.
$$
Furthermore, since $\int_\RR \chi_\delta=1$, we have $\int_\RR \varphi=1$
and $0\le \varphi \le 1/2$.

We now set
$\psi(t):= 2M \varphi(2Mt)$
and define $g(t):=\frac{\pi}{2}\int_{-\infty}^t \psi(s)ds$.
Evidently,
$g\in C^\infty(\RR)$, $\supp g' \subset [-\frac{1}{2}, \frac{1}{2}]$,
$g(t)+g(-t)=\frac{\pi}{2}$ for
$t\in\RR$, $0\le g\le \pi/2$, $\|g'\|_\infty \le (\pi/2)\|\psi\|_\infty \le (\pi/2)M$ and
\begin{equation}\label{est-gk}
\|g^{(k)}\|_\infty \le \frac{\pi}{2}\|\psi^{(k-1)}\|_\infty
\le \frac{\pi}{2}(2M)^{k}k!\cL(k-1)^k
\quad\mbox{for}\quad k\ge 2.
\end{equation}

Apparently $\ha(t):=\frac{2}{\pi}g(\frac 3 2-t)$ is an admissible function of type (a)
and $\ha$ belongs to $\cS(1, \cL; 1, 2M)$.
Also $\ha(t)-\ha(2t)$ is an admissible function of type (b) belonging to $\cS(1,\cL; 2, 4M)$.

To construct an admissible function of type $(c)$
we write $\phi(t):= \sin g(t)$, $t\in \RR$.
From above,
$\phi(t)^2+\phi(-t)^2=1$ for $t\in\RR$.
We define
$$
\ha(t):=
\left\{
\begin{array}{lcl}
\phi(2t-\frac{3}{2}) &\mbox{if}& t\in [\frac{1}{2}, 1],\\
\phi(\frac{3}{2}-t) &\mbox{if}& t\in (1, 2],\\
0 &\mbox{if}& \RR\setminus [\frac{1}{2}, 2].
\end{array}
\right.
$$
We claim that $\ha$ is an admissible cutoff function of type (c)
and $\ha \in \cS(1,\cL; 8, 8M)$.
All required conditions on $\ha$ are trivial to verify  but the estimate
\begin{equation}\label{nontrivial}
\frac{1}{k!}\|\ha^{(k)}\|_\infty
\le 8\big(8M\cL(k-1)\big)^k, \quad k\ge 1.
\end{equation}
Let $t_0\in (-\frac{1}{2}, \frac{1}{2})$ and set
$g_k(t):=\sum_{j=0}^k \frac{(t-t_0)^j}{j!}g^{(j)}(t_0)$.
It is easy to see that $\phi^{(k)}(t_0)= [\sin g_k]^{(k)}(t_0)$
and since $\sin g_k(z)$ is an entire function, by the Cauchy formula,
\begin{equation}\label{Cauchy}
\phi^{(k)}(t_0)=\frac{k!}{2\pi i}\int_\cC \frac{\sin g_k(z)}{(z-t_0)^{k+1}}dz,
\end{equation}
where $\cC:=\{z\in \bC:|z-t_0|=r\}$ with $r=\frac{1}{4M\cL(k-1)}$.
By (\ref{est-gk}) we have for $z\in \cC$ and $k\ge 1$
\begin{align*}
|g_k(z)|&\le \frac{\pi}{2}\Big(1+\frac{M}{4M\cL(k-1)}+\sum_{j=2}^k
\frac{(2M)^j j! \cL(j-1)^j}{j![4M\cL(k-1)]^j}\Big)\\
&\le \frac{\pi}{2}\Big(1+\frac{1}{4}+\sum_{j=2}^k\frac{1}{2^j}\Big)=\frac{7\pi}{8}
\end{align*}
and hence
$|\sin g_k(z)|\le (e^{7\pi/8}+e^{-7\pi/8})/2<8$ for $z\in \cC$.
From this and (\ref{Cauchy}) we get
$$
|\phi^{(k)}(t_0)| 
\le 8 k![4M\cL(k-1)]^k,
$$
which implies (\ref{nontrivial}).
$\qed$


\begin{remark}\label{rem:sharp-est}
{\rm
Theorem~\ref{thm:cutoff-d=1} is sharp in the sense that if
$\int_0^\infty\frac{dt}{(t+1)\cL(t)} =\infty$,
then there is no admissible cutoff function $\ha$ belonging to $\cS(1, \cL; \gamma, \tilde\gamma)$
for any $\gamma, \tilde\gamma >0$.
The argument is precisely the same as in \cite[Remark 3.2]{IPX}.
}
\end{remark}

\subsection{Localized kernels induced by Jacobi polynomials}
\label{Jacobi-kernels}

The Jacobi polynomials $P_n^{(\a,\b)}$, $n=0, 1, \dots$,
form an orthogonal basis for the weighted space
$L^2([-1, 1], w_{\a,\b})$ with weight $w_{\a,\b}(t):=(1-t)^\a(1+t)^\b$.
For various technical reasons we shall assume that
$\a, \b \ge -1/2$.
The Jacobi polynomials are traditionally normalized by
$P_n^{(\a,\b)}(1)=\binom{n+\a}{n}$.
It is well known that~\cite[(4.3.3)]{Sz}
$$
\int_{-1}^1 P_n^{(\a,\b)}(t)
P_m^{(\a,\b)}(t)w_{\a,\b}(t)dt  = \delta_{n, m} h_n^{(\a,\b)},
$$
where
\begin{equation}\label{def-hn}
h_n^{(\a,\b)} =
\frac{2^{\a+\b+1}}{(2n+\a+\b+1)}
\frac{\Gamma(n+\a+1)\Gamma(n+\b+1)}{\Gamma(n+1)\Gamma(n+\a+\b+1)}.
\end{equation}
Hence
\begin{equation}\label{JacobiPoly-l2}
\tP_n^{(\a, \b)}=(h_n^{(\a,\b)})^{-1/2}P_n^{(\a, \b)}
\end{equation}
is the $n$th degree Jacobi polynomial normalized in $L^2([-1, 1], w_{\a, \b})$.

We are interested in kernels of the form
\begin{equation}\label{def.L}
L_n^{\a,\b}(x,y)=\sum_{j=0}^\infty \ha \Big(\frac{j}{n}\Big)
     \tP_j^{(\a,\b)}(x) \tP_j^{(\a,\b)}(y),
\end{equation}
for smooth cutoff functions $\ha: [0, \infty) \mapsto \bC$.

In \cite{PX1} (see also \cite{BD}) it was proved that the kernels $L_n^{\a,\b}(x,y)$ decay rapidly
away from the main diagonal in $[-1, 1]^2$ for compactly supported $C^\infty$ cutoff functions
$\ha$ which are constants around $t=0$.
It was also proved in \cite{IPX} that for such cutoff functions
with ``small" derivatives the localization of these kernels is sub-exponential.
Furthermore, it was shown that the behavior
of $\ha$ at $t=0$ plays a critical role for the localization of $L_n^{\a,\b}(x,y)$,
in particular, the fact that $\ha$ is $C^\infty$ and compactly supported does not
guarantee rapid decay of the kernels $L_n^{\a,\b}(x,y)$.

Here we extend that localization result from \cite{PX1} to smooth cutoff functions $\ha$
with multiple zeros of their first derivatives at $t=0$.
To give this result we need the quantities:
$w_{\a,\b}(0; x):= 1$ and
\begin{equation}\label{Jacobi-weight-n}
w_{\a,\b}(n; x)
:= (1-x + n^{-2})^{\a+1/2} (1+x + n^{-2})^{\b+1/2},
\quad n\ge 1.
\end{equation}
We shall also use the distance
$\rho(x, y):=|\arccos x-\arccos y|$ on $[-1, 1]$.


\begin{theorem}\label{thm:Jacobi-localization-0}
Let $\ha\in C^{3k-1}[0, \infty)$ for some integer $k\ge 1$,
$\supp \ha \subset [0, 2]$, and $\ha^{(m)}(0)=0$ for $m=1, 2, \dots, 3k-1$.
Then there exists a constant $c>0$  of the form
$c=c(k,\a,\b)\|\ha^{(3k-1)}\|_\infty$
such that the kernels from $(\ref{def.L})$ satisfy
\begin{equation}\label{Jacobi-bound0}
|L_n^{\a,\b} (x,y)|
\le c\frac{n}{
\sqrt{w_{\a,\b}(n; x)}\sqrt{w_{\a,\b}(n; y)}}\big(1 + n\rho(x, y)\big)^{-k},
\quad x, y\in [-1, 1].
\end{equation}
Consequently, if $\ha$ is an admissible cutoff function, then the above estimate holds
for any $k\ge1$.
\end{theorem}

As in \cite{PX1} estimate (\ref{Jacobi-bound0}) follows by the localization of $L_n^{\a,\b}(x,1)$
given in the next theorem.
Denote
$\Q_n^{\a,\b}(x) := L_n^{\a,\b}(x,1)$.
It is readily seen that (see e.g.  \cite{PX1})
\begin{equation}\label{eq:L-n}
\Q_n^{\a,\b}(x) =  c^\star
\sum_{j=0}^\infty \ha\Big(\frac{j}{n}\Big)
\frac{(2j+\a+\b+1)\Gamma(j+\a+\b+1)}{\Gamma(j+\b+1)}P_j^{(\a,\b)}(x),
\end{equation}
where
$c^\star:=2^{-\a-\b-1}\Gamma(\a+1)^{-1}$.


\begin{theorem} \label{thm:est-Q}
Let $\ha$ be as in Theorem~\ref{thm:Jacobi-localization-0} and $\a \ge \b \ge -1/2$.
Then for any $r\ge 0$
\begin{equation}\label{est.Ln}
\Big|\frac{d^r}{dx^r} \cQ_n^{\a,\b}(\cos \theta)\Big|
\le c\frac{n^{2\a+2r+2}}{(1+n \theta)^{k}},
\quad 0 \le \theta \le \pi.
\end{equation}
Here $c$ is of the form
$c=c(k, r,\a)\|\ha^{(3k-1)}\|_\infty$.
\end{theorem}

\noindent
{\bf Proof.}
We shall proceed quite similarly as in the proof of Theorem~4.2 in \cite{IPX}
and, therefore, we shall use some notation and facts from that proof.

We shall only prove (\ref{est.Ln}) for $r=0$;
then in general (\ref{est.Ln}) follows by using Markov's inequality as in  \cite{IPX}.

%
%

We trivially have (see (4.8) in \cite{IPX})
$|\Q_n^{\a,\b}(\cos\theta)| \le c n^{2 \a+2}$,
which gives (\ref{est.Ln}) ($r=0$) for $0\le \theta\le 1/n$.


The following identity is crucial in estimating $|\Q_n^{\a,\b}(\cos \theta)|$
\cite[(4.5.3)]{Sz}:
\begin{align}\label{key-ident}
&\sum_{\nu=0}^n \frac{(2\nu+\a+k+\b+1)\Gamma(\nu+\a+k+\b+1)}{\Gamma(\nu+\b+1)}
P_\nu^{(\a+k,\b)}(x)\\
&\qquad\qquad= \frac{\Gamma(n+\a+k+1+\b+1)}{\Gamma(n+\b+1)}
P_n^{(\a+k+1,\b)}(x).\notag
\end{align}
%
%
We define
$
A_0(t):= (2t+\a+\b+1)\ha\left(\frac{t}{n}\right)
$
and inductively
\begin{equation}\label{Ak+1}
A_{k+1}(t):= \frac{A_k(t)}{2t+\a+k+\b+1}-\frac{A_k(t+1)}{2t+\a+k+\b+3},
\quad k\ge 0.
\end{equation}
We apply summation by parts $k$ times starting from (\ref{eq:L-n}) and using every time
(\ref{key-ident}) and \eqref{Ak+1} to obtain
\begin{equation}\label{Ln-parts-k}
\Q_n^{\a,\b}(x) = c^\star\sum_{j=0}^\infty
A_k(j)\frac{\Gamma(j+\a+k+\b+1)}{\Gamma(j+\b+1)}P_j^{(\a+k,\b)}(x).
\end{equation}
Observe first that
$
A_1(t)= \ha(\frac{t}{n})-\ha(\frac{t+1}{n})
=\frac{1}{n}\int_0^1\ha'(\frac{t+s}{n})ds
$
and hence
$
A_1^{(m)}(t)=\frac{1}{n^{m+1}}\int_0^1\ha^{(m+1)}(\frac{t+s}{n})ds,
$
which leads to
$$
|A_1^{(m)}(t)| \le \frac{1}{n^{m+1}} \big\|\ha^{(m+1)}\big\|_{L^\infty[\frac{t}{n}, \frac{t+1}{n}]}.
$$
On the other hand, since  $\ha^{(\ell)}(0)=0$ for $\ell=1, 2, \dots, 3k-1$, then by Taylor's theorem
\begin{equation}\label{Taylor}
|\ha^{(m+1)}(z)| \le \frac{z^{2k-1}}{(2k-1)!}\big\|\ha^{(2k+m)}\big\|_{L^\infty[0, z]}
\quad \mbox{whenever $\; m+1\le k$, $z>0$.}
\end{equation}
Therefore,
\begin{equation}\label{A1}
|A_1^{(m)}(t)| \le \frac{1}{n^{m+1}} \Big(\frac{t+1}{n}\Big)^{2k-1}
\big\|\ha^{(2k+m)}\big\|_{L^\infty[0, \frac{t+1}{n}]},
\quad m+1\le k, \; t >0.
\end{equation}
We next estimate $|A_l^{(m)}(t)|$ by induction on $l$.
We claim that
\begin{equation}\label{Al}
|A_l^{(m)}(t)| \le \frac{c}{(t+1)^{m+2l-1}} \Big(\frac{t+l}{n}\Big)^{2k-1}
\max_{2k \le \ell\le 2k+m+l-1}\big\|\ha^{(\ell)}\big\|_{L^\infty[0, \frac{t+l}{n}]}
\end{equation}
if  $m+l \le k$, $m\ge 0$, $l\ge 1$, and $0\le t \le 2n$, where $c=c(l, m)$,
and hence, using (\ref{Taylor}),
\begin{align}\label{Ak}
|A_k(t)|
&\le \frac{c(k)}{(t+1)^{2k-1}} \Big(\frac{t+k}{n}\Big)^{2k-1}
\max_{2k \le \ell \le 3k-1}\big\|\ha^{(\ell)}\big\|_{L^\infty[0, \frac{t+k}{n}]}\\
&\le \frac{c(k)}{n^{2k-1}} \big\|\ha^{(3k-1)}\big\|_\infty.\notag
\end{align}
Indeed, estimate (\ref{A1}) gives (\ref{Al}) for $l=1$.
Suppose  (\ref{Al}) holds for some $l\ge 1$ and all $m\ge 0$ such that  $m+l \le k$.
Then by (\ref{Ak+1})
$
A_{l+1}(t)= -\int_0^1 G_l'(t+s)ds
$
with $G_l(t):= \frac{A_l(t)}{2t+\a+l+\b+1}$
and hence
$
A_{l+1}^{(m)}(t)= -\int_0^1 G_l^{(m+1)}(t+s)ds.
$
We have
$$
G_l^{(m+1)}(t)
=\sum_{\nu=0}^{m+1}\binom{m+1}{\nu}A_l^{(\nu)}(t)
\frac{(-2)^{m+1-\nu}(m+1-\nu)!}{(2t+\a+l+\b+1)^{m+2-\nu}}
$$
and using the inductive assumption
\begin{align*}
|A_{l+1}^{(m)}(t)|
& \le c \max_{2k \le \ell \le 2k+m+l}\big\|\ha^{(\ell)}\big\|_{L^\infty[0, \frac{t+l+1}{n}]}\\
&\qquad\qquad\qquad\qquad
\times\sum_{\nu=0}^{m+1} \frac{1}{(t+1)^{\nu+2l-1}}
\Big(\frac{t+l+1}{n}\Big)^{2k-1}\frac{1}{(t+1)^{m+2-\nu}}\\
&\le \frac{c}{(t+1)^{m+2l+1}} \Big(\frac{t+l+1}{n}\Big)^{2k-1}
\max_{2k \le \ell \le 2k+m+l}\big\|\ha^{(\ell)}\big\|_{L^\infty[0, \frac{t+l+1}{n}]},
\end{align*}
which confirms (\ref{Al}).


We next prove (\ref{est.Ln}) ($r=0$) for $1/n \le \theta\le \pi/2$.
By (\ref{def-hn}) it readily follows that
$h_n^{(\a+k, \b)}\le c2^k/n$
and it is well known that  (see e.g. (4.18) in \cite{IPX})
$$
|P_n^{(\a+k,\b)}(\cos \theta)| \le \frac{c}{n^{1/2}\theta^{k+\a+1/2}},
\quad 0<\theta\le \pi/2.
$$
We use the above and (\ref{Ak}) in (\ref{Ln-parts-k}) to obtain
for $1/n \le \theta \le \pi/2$
\begin{align*}
|\Q_n^{\a,\b}(\cos \theta)|
\le c \|\ha^{(3k-1)}\|_\infty
\sum_{j=1}^{2n}\frac{ j^{\a+k} }{n^{2k-1} j^{1/2}\theta^{k+\a+1/2}}
\le c \|\ha^{(3k-1)}\|_\infty\frac{n^{2\a+2}}{(n\theta)^{k+\a+1/2}}.
\end{align*}
Hence, estimate (\ref{est.Ln}) (with $r=0$) holds for $1/n\le \theta\le \pi/2$.


Let $\pi/2 < \theta\le \pi-1/n$.
Similarly as in \cite{IPX}
$$
|P_n^{(\a+k, \b)}(\cos\theta)|
\le c2^kn^\beta,
\quad \pi/2\le \theta\le \pi-1/n.
$$
Combining this with (\ref{Ln-parts-k}) and (\ref{Ak}) we get for $\pi/2\le \theta\le \pi-1/n$
\begin{align*}
|\Q_n^{\a,\b}(\cos \theta)|
\le  c  \|\ha^{(3k-1)}\|_\infty n^{-2k+1} \sum_{j=1}^{2n} j^{\b+\a+k}
\le  c  \|\ha^{(3k-1)}\|_\infty \frac{n^{\a+\b+2}}{n^k},
\end{align*}
which implies (\ref{est.Ln}).

In the case $\pi-1/n \le \theta\le \pi$ estimate  (\ref{est.Ln}) follows from the above estimate
exactly as in \cite{IPX}.
This completes the proof of estimate  (\ref{est.Ln}) in the case $r=0$.
$\qed$

\medskip

Estimate (\ref{Jacobi-bound0}) can be improved for admissible cutoff functions which
are constant around $t=0$ and have ``small" derivatives as in Theorem~\ref{thm:cutoff-d=1}:


\begin{theorem}\label{thm:Jacobi-localization-1}
Let $\cL$ and  $M$ be as in \eqref{def-M}.
Suppose $\ha$ is an admissible cutoff function of type $(a)$,  $(b)$ or  $(c)$
which belongs to $\cS(1,\cL; \gamma,\bargamma M)$ for some $\gamma,\bargamma>0$
$($see Theorem~\ref{thm:cutoff-d=1}$)$.
Then the kernels from $(\ref{def.L})$ satisfy
\begin{equation}\label{Jacobi-bound1}
|L_n^{\a,\b} (x,y)|
\le \frac{cn}{
\sqrt{w_{\a,\b}(n; x)}\sqrt{w_{\a,\b}(n; y)}}
\exp\Big\{- \frac{\tilde{c}n\rho(x, y)}{\cL(n\rho(x, y))}\Big\}
\end{equation}
for $x, y\in [-1, 1]$,
where $\tilde{c}=c'/\bargamma M$ with $c'>0$ being an absolute constant
and $c$ depends continuously only on $\a$, $\b$, $\gamma,\bargamma$ and $M$.

In particular, the above result holds for
$\cL_{\ell,\eps}(t)$ from \eqref{log1} with $M=c(\ell)\eps^{-1}$.
\end{theorem}

For the proof of this theorem one first uses Theorem~\ref{thm:cutoff-d=1} to prove
the following estimate for the kernels from \eqref{eq:L-n} with $\ha$ from above
\begin{equation*}
\Big|\frac{d^r}{dx^r} \Q_n^{\a,\b}(\cos \theta)\Big|
\le c n^{2\a+2r+2}\exp\left\{-\frac{\tilde{c}n\theta}{\cL(n\theta)}\right\},
\quad 0 \le \theta \le \pi,
\end{equation*}
and then proceeds exactly as in the proof of Theorem~4.1 in \cite{IPX}.
The proofs are nearly identical to the ones in \cite{IPX} and will be omitted.


\begin{remark}\label{rem:Jacobi-localization-1}
{\rm Theorem \ref{thm:Jacobi-localization-1} remains true if we require that
$\ha\in \cS(1,\cL; \gamma, \bargamma M)$, $\supp \ha \subset [0, 2]$, and
$\ha$ be a constant on $[0,\delta]$ for a fixed $\delta\in(0,1)$.
Then $c$ and $c'$ will depend on $\delta$ as well.
However, we are not aware if Theorem \ref{thm:Jacobi-localization-1}  holds for admissible
cutoff functions which are not constants around $t=0$.
The method of proof of Theorem~\ref{thm:est-Q}
does not give such a result for all admissible cutoff functions.

In \cite{IPX} results similar to Theorem~\ref{thm:Jacobi-localization-1} are proved on
the sphere, ball, simplex and in the context of Hermite and Jaguerre functions with
$\cL$ replaced by $\cL_{\ell, \eps}$. We would like to point out here that with the same proofs
these results hold for a general function $\cL$ as above.
}
\end{remark}

We shall need


\begin{lemma}\label{lem:simple-est}
There exists a constant $c$ depending only on $\a, \b$ such that
\begin{equation}\label{est.Pn2}
|\tP_n^{(\a,\b)}(x)|\le \frac{c}{\sqrt{w_{\a,\b}(n; x)}},\quad\quad x\in[-1, 1],~~n\ge 1.
\end{equation}
\end{lemma}

\noindent
{\bf Proof.}
For $x\in[-1+n^{-2},1-n^{-2}]$ using
\begin{equation*}
    w_{\a,\b}(n; x)\le 2^{\a+\b+1}(1-x)^{\a+1/2}(1+x)^{\b+1/2}
\end{equation*}
we get \eqref{est.Pn2} from the inequality
\begin{equation*}
\sup_{x\in[-1, 1]}(1-x)^{\a+1/2}(1+x)^{\b+1/2}|\tP_n^{(\a,\b)}(x)|^2
\le \frac{2e}{\pi}\big(2+\sqrt{\a^2+\b^2}\big)
\end{equation*}
established in \cite[Theorem 1]{EMN}.
For the remaining $x$ estimate (\ref{est.Pn2}) follows from above invoking
Theorem 8.4.8 in \cite[p. 108]{Di-To}.
$\qed$

The next theorem shows that the kernels $L_n^{\a,\b}$ from (\ref{def.L})
are Lip 1 with respect to the distance $\rho(\cdot, \cdot)$.


\begin{theorem}\label{thm:Jacobi-Lip-1}
Let $\ha\in C^{3k-1}[0, \infty)$ for some
$k > 2\a+2\b+5$,
$\supp \ha \subset [0, 2]$, and $\ha^{(r)}(0)=0$ for $r=1, 2, \dots, 3k-1$.
Then for any $x, y, \xi\in [-1, 1]$ such that $\rho(x, \xi) \le  c_*n^{-1}$,
$n\ge 1$, $c_*>0$, the kernel $L_n^{\a,\b}$ from $(\ref{def.L})$ satisfies
\begin{equation}\label{Jacobi-Lip-1}
|L_n^{\a,\b} (x,y)-L_n^{\a,\b} (\xi,y)|
\le \frac{cn^2\rho(x, \xi)}
{\sqrt{w_{\a,\b}(n; x)}\sqrt{w_{\a,\b}(n; y)}}\big(1 + n\rho(x, y)\big)^{-\sigma},
\end{equation}
where $\sigma=k-2\a-2\b-5$ and $c$ depends only on $k$, $\a$, $\b$, $c_*$, and
$\|\ha^{(3k-1)}\|_\infty$.
Consequently, if $\ha$ is an admissible cutoff function, then the above estimate holds
for any $\sigma > 0$.
\end{theorem}

The proof of this theorem for $\alpha, \beta > -1/2$ utilizes estimate \eqref{est.Ln}
and is identical with the proof of Theorem 2.2 in \cite{KPX1}.
The limit cases $\a=-1/2$ or $\b=-1/2$ are treated as in the proof of \cite[Theorem 4.1]{IPX}.
We omit the details.

\section{Multivariate cutoff functions}\label{multi-kernels}
\setcounter{equation}{0}

As was explained in the introduction, cutoff functions in d-variables will play
a prominent role in the development of weighted $F$- and $B$-spaces on $[-1, 1]^d$.
In this section we introduce two kinds of admissible d-dimensional cutoff functions
and give several constructions of such functions.

\subsection{Admissible d-dimensional cutoff functions}
\label{d-dimensional cutoff}
\setcounter{equation}{0}

To define multivariate cutoff functions we need to introduce some convenient notation.
Given $1\le k\le d$ we define
$\proj_k :\RR^d\to\RR^{d}$ by
\begin{equation}\label{def-proj-k}
\proj_k(t_1,\dots,t_d):=(t_1,\dots,t_{k-1},0,t_{k+1},\dots,t_d).
\end{equation}
We also denote by $\CB_p$ the part of the unit ball of the standard $\ell^p(\RR^d)$ norm
contained in the first octant, i.e.
\begin{equation*}
\CB_p:=\{t\in [0, \infty)^d:~\|t\|_p\le 1\},\quad\quad 1\le p\le\infty.
\end{equation*}


\begin{definition}\label{cutoff-first}
A cutoff function $\hA\in C^\infty [0, \infty)^d$ is said to be {\bf admissible of first kind}
or simply {\bf admissible}
if $\supp \hA \subset [0, 2]^d$ and $\hA$ obeys the First
Boundary Condition, introduced in \S\ref{remedy}, i.e.
for any $t\in [0, \infty)^d$ of the form
$t=\proj_k t$ for some $1\le k\le d$ we have
$D_k^m \hA(t)=0$ for $m=1, 2, \dots$.

Furthermore, $\hA$ is said to be of type $(a)$, $(b)$, or $(c)$ if in addition

\smallskip
$(a)$ $\hA(t)=1$ for $t\in \CB_1$,

\smallskip

$(b)$ $\hA(t)=0$ for $t\in\frac{1}{2} \CB_1$, or

\smallskip
$(c)$ $\hA$ is of type $(b)$ and $\;\sum_{j=0}^\infty |\hA(2^{-j}t)|^2=1\;$
for $\;t\in [0, \infty)^d \setminus [0, 1)^d$.

\end{definition}


\begin{definition}\label{cutoff-second}
A cutoff function $\hA\in C^\infty [0, \infty)^d$ is said to be {\bf admissible of second kind}
and type $(a)$, $(b)$, or $(c)$ if $\supp\hA\subset [0, 2]^d$, $\hA(t)=\hA(\proj_k t)$
for every $t=(t_1,\dots,t_d)\in[0, \infty)^d$ such that $t_k\le\frac{1}{2}\|t\|_\infty$ and

\smallskip
$(a)$ $\hA(t)=1$ if $t\in \CB_1$,

\smallskip

$(b)$ $\hA(t)=0$ if $t\in\frac{1}{2} \CB_1$, or

\smallskip
$(c)$ $\hA$ is of type $(b)$ and $\;\sum_{j=0}^\infty |\hA(2^{-j}t)|^2=1\;$
for $\;t\in [0, \infty)^d \setminus [0, 1)^d$.
\end{definition}

We first note that if $\hA$ is admissible of second kind, then for all $\tau\in[0, \infty)^d$
such that $\tau=\proj_k \tau$ for some $1\le k \le d$, the function $\hA(t)$
is a constant on the segment $t\in[\tau,\tau+\frac{1}{4d-2}e_k]$,
i.e. $\hA$ obeys the Second Boundary Condition from
\S\ref{remedy} with constant $c_*=\frac{1}{4d-2}$.
Consequently, any admissible cutoff function of second kind is admissible of first kind
as well.
To see the above one simply has to observe that $\frac{1}{4d-2}$ is the $\ell^\infty$-distance
of the set
\begin{equation*}
\Big(\bigcap_{k=1}^d \Big\{t_k\ge\frac{1}{2}\|t\|_\infty\Big\}\Big)\bigcap\Big\{\|t\|_1
\ge \frac{1}{2}\Big\}
\end{equation*}
from the coordinate hyperplanes. Note also that for $d=1$ the set of admissible cutoff functions
of first kind and type (a) coincides with the set of admissible cutoff functions of second kind
and type (a); the same for types (b) and (c).


\begin{remark}\label{raeson}
{\rm
As was explained in the introduction the fact that the admissible cutoff functions
satisfy the First Boundary Condition (see \S\ref{remedy}) is crucial
for the rapid decay of the associated tensor product Jacobi polynomial kernel from \eqref{def.L};
this will be established in the next section.

An important reason for introducing admissible cutoff functions of second kind
is that such cutoff functions with ``small" derivatives (\S\ref{small1}, \S\ref{small2})
allow to construct tensor product Jacobi polynomial kernel of sub-exponential localization
(see Theorems~\ref{thm:product-Jacobi-1} and \ref{thm:product-ball-Jacobi-1}),
while as for now we are unable to achieve such localization
with admissible cutoff functions of first kind.
}
\end{remark}


It is easy to construct admissible cutoff functions of type~(a) as products of univatiate
cutoff functions of type~(a).


\begin{lemma}\label{lem:product_a}
Let $\ha_j$, $j=1,\dots,d$, be admissible univariate functions of type $(a)$.
Then
\begin{equation}\label{cutoff-d3}
    \hA(t)=\prod_{j=1}^d \ha_j(t_j)
\end{equation}
is an admissible d-dimensional cutoff function of second kind and type $(a)$.
\end{lemma}

\noindent
{\bf Proof.}
By the definition evidently $\hA\in C^\infty[0,\infty)^d$,  $\supp\hA\subset 2\CB_\infty$
and $\hA(t)=1$ if $t\in \CB_\infty\supset\CB_1$.
Furthermore, $\hA(t)=\hA(\proj_k t)$ for all $t\in [0, \infty)^d$
such that $\|t\|_\infty< 2$ and $t_k\le 1$.
From this and $\supp\hA\subset 2\CB_\infty$ it follows that  $\hA(t)=\hA(\proj_k t)$
for all $t=(t_1,\dots,t_d)\in[0, \infty)^d$ such that $t_k\le\frac{1}{2}\|t\|_\infty$.
$\qed$

\smallskip

The construction of admissible cutoff functions of type~(b) is straightforward using
admissible cutoff functions of type~(a):


\begin{lemma}\label{lem0}
If $\hA_1$, $\hA_2$ are admissible of type $(a)$ $($any kind$)$, then
\begin{equation}\label{cutoff-d4}
    \hA(t)=\hA_1(t)-\hA_2(2t)
\end{equation}
is admissible of type $(b)$.
\end{lemma}

For the definition of $F$- and $B$-spaces on $[-1, 1]^d$ we shall utilize
admissible cutoff functions $\hA$ of type (b) with the property that
the dyadic dilations of $\supp \hA$ essentially cover the whole spectrum.
More precisely, we shall need admissible functions $\hA$, which obey the following
{\bf dyadic covering condition}:
\begin{equation}
\begin{aligned}\label{star}
&\mbox{For any $t\in [0,\infty)^d$ with $\|t\|_\infty=1$ there is $0< \gamma \le 1$ such that}\\
& \hspace{1.2in} \inf_{\lambda\in[\gamma,2\gamma]}|\hA(\lambda t)|\ge c>0.
\end{aligned}
\end{equation}

\noindent
Note that this condition yields
$\sum_{j=0}^\infty |\hA(2^{-j}t)| \ge c>0$ for $t\in [0, \infty)^d\setminus \CB_\infty$
which justifies our terminology.

From the constructions of admissible cutoff functions below it will be clear that
it is easy to construct admissible functions $\hA$ of type (b) which satisfy
condition~(\ref{star}).

The following lemma will be instrumental in the development of $F$- and $B$-spaces.


\begin{lemma}\label{lem:dual-cutoff}
For any admissible function $\hA$ of first or second kind and type $(b)$ satisfying
the dyadic covering condition $(\ref{star})$
there exists an admissible function $\hB$ of type $(b)$ $($and the same kind$)$ such that
\begin{equation}\label{conjugate1}
\sum_{j=0}^\infty \overline{\hA(2^{-j}t)}\hB(2^{-j}t)=1
\quad \mbox{for $\;t\in [0, \infty)^d\backslash \CB_\infty$.}
\end{equation}
\end{lemma}

\noindent
{\bf Proof.} We shall only prove this lemma for an admissible function of second kind,
since the case of first kind cutoff functions is easier.

We define $\hB(t):=0$ for $t\in\frac{1}{2}\CB_1$ and
$t\in [0, \infty)^d\backslash 2\CB_\infty$.
For the remaining $t\in [0, \infty)^d$ we set
\begin{equation}\label{conjugate2}
\hB(t):=\frac{\hA(t)}{\sum_{j=-\infty}^\infty |\hA(2^{-j}t)|^2}.
\end{equation}
For every $t\in [0, \infty)^d$ the sum in the denominator of \eqref{conjugate2} is non-zero
on account of property (\ref{star}) and contains no more that $2+\log_2 d$ non-zero terms.
Hence $\hB \in C^\infty[0, \infty)^d$.
On the other hand, for $t\in [0, \infty)^d\backslash \CB_\infty$
we have $2^{-j}t\notin 2\CB_\infty$ for $j<0$ and
the sum in the denominator of \eqref{conjugate2} reduces to $j\ge 0$.
Hence \eqref{conjugate1} is trivially satisfied.
Finally, if $t\in[0, \infty)^d$ and
$t_k\le\frac{1}{2}\|t\|_\infty$ for some $1\le k \le d$,
then $2^{-j}t_k\le\frac{1}{2}\|2^{-j}t\|_\infty$
for $j\in \ZZ$ and $\hA(2^{-j}t)=\hA(\proj_k (2^{-j}t))$, which implies $\hB(t)=\hB(\proj_k t)$.
\qed

\smallskip

The construction of admissible cutoff functions of type (c) will require some care.
We shall give several constructions of cutoff functions below.

\subsection{Construction of admissible cutoff functions via quasi-norms}

One approach for constructing admissible d-dimensional cutoff functions is based on the following lemma.


\begin{lemma}\label{lem3}
Suppose the function $\cN:\RR^d\to\RR$ is in $C^\infty(\RR^d\backslash\{0\})$
and for $t\in\RR^d$ obeys
\begin{gather}
\cN(\alpha t)=\alpha \cN(t),\quad\alpha>0,\label{cond1}\\
\|t\|_\infty\le \cN(t)\le\|t\|_1,\label{cond2}\\
\cN(t)=\cN(\proj_m t) \quad\mbox{provided}\quad\ |t_m|\le\frac{1}{2}\|t\|_\infty,\quad
m=1,\dots,d.\label{cond3}
\end{gather}
If $\ha$ is an admissible univariate function of type $(a)$, $(b)$ or $(c)$, then
\begin{equation*}
    \hA(t)=\ha(\cN(t))
\end{equation*}
is an admissible d-dimensional function of second kind and type $(a)$, $(b)$ or $(c)$, respectively.
\end{lemma}

The proof is straightforward.

A simple way to construct a function $\cN$ satisfying the conditions of Lemma \ref{lem3}
is the following. Let $\hc$ be an even real-valued function, whose restriction on $[0,\infty)$
is an admissible univariate function of type $(a)$, satisfying $0\le \hc\le 1$.
For $d\in\NN$ and $t\in\RR^d\backslash\{0\}$ set
\begin{equation}\label{quasin1}
\cN(t):=\sum_{m=1}^d |t_m| \prod_{j=1}^d \hc\Big(\frac{t_j}{t_m}\Big),
\end{equation}
where $\hc\big(\frac{\tau}{0}\big):=0$ for every real $\tau$, including $\tau=0$.
For $t=0$ by continuity we set $\cN(0)=0$.

Given $t\ne 0$ let $k$ be such that
$\|t\|_\infty = |t_k|$. If $|t_m|\le\frac{1}{2}\|t\|_\infty$,
then $|t_k/t_m|\ge 2$ and hence $\hc\big(\frac{t_k}{t_m}\big)=0$.
Observing also that $\hc\big(\frac{t_j}{t_m}\big)=1$ if $|t_j|\le|t_m|$
we see that \eqref{quasin1} can be rewritten as
\begin{equation}\label{quasin2}
\cN(t)=\sum_{1\le m \le d, \; |t_m|>\frac{1}{2}\|t\|_\infty} |t_m|
\prod_{1\le j\le d,\; |t_j|>|t_m|} \hc\Big(\frac{t_j}{t_m}\Big).
\end{equation}
It follows immediately from \eqref{quasin2} that $\cN$ belongs to
$C^\infty(\RR^d\backslash\{0\})$ and satisfies condition \eqref{cond1}.
If $m$ is such that  $|t_m|\le\frac{1}{2}\|t\|_\infty$,
then $|t_m|$ does not participate in the right-hand side of \eqref{quasin2}
and hence $\cN$ satisfies \eqref{cond3}.
The inequality $\cN(t)\le\|t\|_1$ follows from $0\le \hc\le 1$.
Finally, from \eqref{quasin2} and $0\le \hc\le 1$ we get $\cN(t)\ge |t_k|=\|t\|_\infty$
and thus \eqref{cond2} is also satisfied.
Thus, we have proved


\begin{corollary}\label{cor-admis-3d}
Let $\cN$ be given by \eqref{quasin1}, where $\hc$ is an even real-valued function,
whose restriction on $[0,\infty)$ is an admissible univariate function of type $(a)$,
satisfying $0\le \hc\le 1$. If $\ha$ is an admissible univariate function
of type $(a)$, $(b)$ or $(c)$, then
\begin{equation*}
\hA(t)=\ha(\cN(t))
\end{equation*}
is an admissible d-dimensional function of second kind and type $(a)$, $(b)$ or $(c)$, respectively.
\end{corollary}


\subsection{Construction of admissible d-dimensional cutoff functions via norms}

From \eqref{quasin1}-\eqref{quasin2} it follows that $\cN$ is a quasi-norm.
A necessary and sufficient condition for $\cN$ to be a norm is the convexity of
the unit ball $\CB=\{t:~\cN(t)\le 1\}$.

The construction of the boundary $\partial\CB=\{t:~\cN(t)= 1\}$ of the unit ball of
a norm $\cN$ satisfying the conditions of Lemma \ref{lem3} can be carried out by induction
on the dimension.
First, one gets the boundaries of the $d-1$ dimensional unit balls on every coordinate hyperplane.
Second, one extends them into the first octant by line segments of length $\frac{1}{2}$.
Third, one completes the surface of the unit ball boundary in the first octant
by convex $C^\infty$ blending.
Finally, one extends it by symmetry to the remaining octants and defines the norm
from the ball in a standard way.

If instead of  convex $C^\infty$ blending in the above scheme it is used a $C^\infty$
blending satisfying \eqref{cond2}, then one obtains a quasi-norm $\cN$ satisfying
all conditions of Lemma~\ref{lem3}.
We shall not further elaborate on this construction.

\subsection{Construction of admissible cutoff functions via quasi-norms with ``small'' derivatives}
\label{small1}
In analogy to Theorem \ref{thm:cutoff-d=1} we construct here admissible d-dimensional cutoff functions
with ``small'' derivatives.

In this construction we shall utilize classes of $C^\infty$ functions of this type:
\begin{multline*}
\cR(a,b,\PsiZ;\gamma,\bargamma)
:=\Big\{f\in C^\infty[a,b]:~
\frac{1}{k!}\|f^{(k)}\|_{L^\infty[a,  b]}\le \gamma\big(\bargamma \PsiZ(k)\big)^k,\quad \forall k\in\NN\Big\},
\end{multline*}
where $\PsiZ$ is a given positive non-decreasing function defined at least on $\NN$, and
$\gamma,\bargamma>0$ are parameters independent of $k$.
Obviously, the sum and the product of two functions from such classes also belong to a class
like that (as the parameters $\gamma,\bargamma$ may vary).
More importantly, the composition of two functions also belongs to such a class as
the following lemma shows.


\begin{lemma}\label{lem:CompositeDerivative1}
Let $\PsiZ(v)>0$ for $v\in[1,\infty)$ and let $\ln\PsiZ(v)^v$ be convex on $[1, \infty)$.
If $f\in\cR(a_1,b_1,\PsiZ;\gamma_1,\bargamma_1)$,
$g\in\cR(a_2,b_2,\PsiZ;\gamma_2,\bargamma_2)$ and the range of $f$ is in $[a_2,b_2]$,
then the composition $g\circ f\in\cR(a_1,b_1,\PsiZ;\gamma_2,\bargamma_1(\gamma_1\bargamma_2\PsiZ(1)+1))$.
\end{lemma}

\noindent
{\bf Proof.}
In order to find an estimate for $D^k(g\circ f)$ we apply Fa\`{a} di Bruno's formula in the form
\begin{equation}\label{eq:CompositeDerivative1}
\frac{1}{k!}D^k(g\circ f)
=\sum_{m\in\EE_k}\frac{|m|!}{m_1!\dots m_k!}\frac{(D^{|m|}g)\circ f}{|m|!}
\prod_{j=1}^k\Big(\frac{D^j f}{j!}\Big)^{m_j},
\end{equation}
where
\begin{equation*}
\EE_k=\{m\in\NN_0^k:~\sum_{j=1}^k j m_j=k\}.
\end{equation*}
Note that $\sum_{j=1}^k j m_j=k$ implies that at most $O(\sqrt{k})$ of $m_j$ can be non-zero.
We assume that $(\cdot)^0=1$ in the product
in \eqref{eq:CompositeDerivative1} even if the argument may be $0$.

Applying the estimates on the derivatives of $f$ and $g$ we get from \eqref{eq:CompositeDerivative1}
\begin{align}\label{eq:CompositeDerivative2}
&\frac{1}{k!}\|D^k(g\circ f)\|_\infty \\
&\le \sum_{m\in\EE_k}\frac{|m|!}{m_1!\dots m_k!}\gamma_2(\bargamma_2 \PsiZ(|m|))^{|m|}
\prod_{j=1}^k \gamma_1^{m_j}(\bargamma_1 \PsiZ(j))^{jm_j}\nonumber\\
&= \gamma_2\bargamma_1^k \sum_{n=1}^k (\gamma_1\bargamma_2)^n \PsiZ(n)^n
\sum_{m\in\EE_k,~|m|=n} \frac{n!}{m_1!\dots m_k!} \prod_{j=1}^k \PsiZ(j)^{jm_j}\nonumber.
\end{align}
Now, from the convexity of $v\ln\PsiZ(v)$ we get
\begin{equation}\label{eq:CompositeDerivative3}
\ell\ln\PsiZ(\ell)+j\ln\PsiZ(j)\le \ln\PsiZ(1)+(\ell+j-1)\ln\PsiZ(\ell+j-1)
\quad \forall \ell, j>1.
\end{equation}
If in a multi-index $m$ we increase $m_1$ and $m_{\ell+j-1}$ by $1$
and decrease $m_\ell$ and $m_j$ by $1$, then the quantities
$\sum_{j=1}^k j m_j$ and $\sum_{j=1}^k m_j$ remain unchanged.
Observing that this operation decreases $\sum_{j=2}^k m_j$ by $1$
and applying inductively \eqref{eq:CompositeDerivative3} we obtain that
among all $m\in\EE_k$ with $|m|=n$ the largest value of the product
$\prod_{j=1}^k \PsiZ(j)^{jm_j}$ is attained for $m_1=n-1, m_{k-n+1}=1$,
and $m_j=0$ if $j\ne 1$ and $j\ne k-n+1$,
i.e.
\begin{equation}\label{eq:CompositeDerivative4}
    \prod_{j=1}^k \PsiZ(j)^{jm_j}\le \PsiZ(1)^{n-1}\PsiZ(k-n+1)^{k-n+1}.
\end{equation}
Using \eqref{eq:CompositeDerivative4} and
\begin{equation*}
\sum_{m\in\EE_k,\; |m|=n} \frac{n!}{m_1!\dots m_k!}=\binom{k-1}{n-1}
\end{equation*}
(see e.g. \cite[Section 5.5]{Rio}) in \eqref{eq:CompositeDerivative2}
and further applying \eqref{eq:CompositeDerivative3} with $\ell=n$, $j=k-n+1$
we finally get
\begin{align*}
\frac{1}{k!}\|D^k(g\circ f)\|_\infty
&\le \gamma_2\bargamma_1^k \sum_{n=1}^k (\gamma_1\bargamma_2)^n \PsiZ(n)^n
\binom{k-1}{n-1} \PsiZ(1)^{n-1}\PsiZ(k-n+1)^{k-n+1}\nonumber\\
&\le \gamma_2\bargamma_1^k \PsiZ(k)^k \sum_{n=1}^k  \binom{k-1}{n-1}
\PsiZ(1)^{n} (\gamma_1\bargamma_2)^n \nonumber\\
&= \gamma_2\bargamma_1^k \PsiZ(k)^k \gamma_1\bargamma_2
\PsiZ(1)(\gamma_1\bargamma_2\PsiZ(1)+1)^{k-1} \nonumber\\
&\le \gamma_2\big[\bargamma_1(\gamma_1\bargamma_2\PsiZ(1)+1) \PsiZ(k)\big]^k  \nonumber.\qed
\end{align*}


We shall utilize Lemma~\ref{lem:CompositeDerivative1} to the composition of admissible cutoff
functions with ``small" derivatives in the sense of Definition~\ref{Def-S}
(see Theorem~\ref{thm:cutoff-d=1}), where $\cL$ obeys an additional {\bf convexity condition}.
Namely, we shall assume that
\begin{equation}\label{def-M-conv}
\cL~\mbox{satisfies \eqref{def-M} and}~(t+1)\ln\cL(t)~\mbox{is convex on}~[0,\infty).
\end{equation}

The functions $\cL_{0, \eps}$ and $\cL_{\ell, \eps}$ from \eqref{log1} are examples of functions
$\cL$ satisfying this condition. 


\begin{theorem}\label{thm:cutoff2}
Let $\cL$ satisfy \eqref{def-M-conv} and let $M$ be given by \eqref{def-M}.
Then the set $\cS(d,\cL;8,10d(2d-1)M(8(d+2)M+1))$ contains an admissible cutoff function
$\hA$, $0\le \hA\le 1$, of second kind and any type:
$(a)$, $(b)$ or $(c)$ $($see Definition~\ref{cutoff-second}$)$.
\end{theorem}

\noindent
{\bf Proof.}
Set $\hA(t)=\ha(\cN(t))$, where $\ha\in\cS(1,\cL;8,8M)$, $0\le\ha\le1$,
is an admissible function of any type $(a)$, $(b)$ or $(c)$ from Theorem~\ref{thm:cutoff-d=1}
 and $\cN=\cN_d$ is given by \eqref{quasin1}
with $\hc|_{[0,\infty)}\in\cS(1,\cL;1,2M)$, $0\le\hc\le1$,
being an admissible function of type $(a)$ from Theorem~\ref{thm:cutoff-d=1}.
Then $\hA$ is an admissible multivariate cutoff function of the same type as $\ha$ according to
Corollary \ref{cor-admis-3d}. Moreover, $0\le \hA\le 1$.

In estimating $D_j^k \hA(t_1,\dots,t_d)$ we may assume without loss of generality that $j=d$.
Further, we consider only
\begin{equation}\label{bound_t}
    \frac{1}{4d-2}\le t_d\le 2,
\end{equation}
because $D_d^k \hA(t_1,\dots,t_d)=0$ if $0\le t_d < \frac{1}{4d-2}$ or $t_d>2$.

In order to apply Lemma \ref{lem:CompositeDerivative1} with $g=\ha$ and $f=\cN$
(as a function of $t_d$) we need upper bounds for $D_d^k \cN(t_1,\dots,t_d)$.
From \eqref{quasin1} we write
\begin{equation*}
    \cN(t)=:\sum_{m=1}^d F_m(t),\quad F_m(t):=t_m \prod_{j=1}^d \hc\Big(\frac{t_j}{t_m}\Big).
\end{equation*}
 Then for $m=1,\dots,d-1$ we have
\begin{equation*}
    D_d^k F_m(t)=(D^k\hc)\Big(\frac{t_d}{t_m}\Big)t_m^{1-k} \prod_{j=1}^{d-1} \hc\Big(\frac{t_j}{t_m}\Big),
\end{equation*}
which on account of \eqref{bound_t}
and since $D^k\hc(\tau)= 0$ for $\tau\notin[1,2]$ and $0\le\hc\le 1$ implies
\begin{equation}\label{bound_m}
\frac{1}{k!}|D_d^k F_m(t)|
\le \big(8(2d-1)M\cL(k-1)\big)^k
\end{equation}
for all $(t_1,\dots,t_{d-1})\in[0,\infty)^{d-1}$.
Using the formulas for derivatives of a product  we get
\begin{align}\label{bound_d}
D_d^k F_d(t)&=t_d D_d^{k} \Big(\prod_{j=1}^{d-1} \hc\Big(\frac{t_j}{t_d}\Big)\Big)
+k D_d^{k-1} \Big(\prod_{j=1}^{d-1} \hc\Big(\frac{t_j}{t_d}\Big)\Big)\\
&=t_d \sum_{|m|=k} \frac{k!}{m_1!\dots m_{d-1}!}\prod_{j=1}^{d-1} D_d^{m_j}
\Big(\hc\Big(\frac{t_j}{t_d}\Big)\Big)\nonumber\\
&+ \sum_{|m|=k-1} \frac{k!}{m_1!\dots m_{d-1}!}
\prod_{j=1}^{d-1} D_d^{m_j}\Big(\hc\Big(\frac{t_j}{t_d}\Big)\Big).\nonumber
\end{align}
For $1\le t_j/t_d\le 2$ the function $t_j/t_d$ of $t_d$ belongs to $\cS(1,\cL;2,2(2d-1))$.
Hence, by Lemma \ref{lem:CompositeDerivative1} with $\PsiZ(v)=\cL(v-1)$ it follows that
$\hc(t_j/t_d)\in \cS(1,\cL;1,10(2d-1)M)$.
Using this and \eqref{bound_t} in \eqref{bound_d} we get
\begin{align}\label{bound_d2}
\frac{1}{k!}|D_d^k F_d(t)|
    &\le 2\sum_{|m|=k} \prod_{j=1}^{d-1}\big(10(2d-1)M
\cL(m_j-1)\big)^{m_j}\\
    &+ \sum_{|m|=k-1}\prod_{j=1}^{d-1}\big(10(2d-1)M
\cL(m_j-1)\big)^{m_j}\nonumber\\
    &\le 3\Big(\sum_{|m|=k} 1 \Big)\big(10(2d-1)M
\cL(k-1)\big)^{k}\nonumber\\
    &\le 3\big(10d(2d-1)M\cL(k-1)\big)^{k}.\nonumber
\end{align}
We recall that the terms in \eqref{bound_d2} with $m_j=0$ are considered equal $1$.
Now, combining \eqref{bound_m} and \eqref{bound_d2} we get for all
$(t_1,\dots,t_{d-1})\in[0,\infty)^{d-1}$ and $t_d$ in \eqref{bound_t}
\begin{equation}\label{bound_N}
    \frac{1}{k!}|D_d^k \cN(t)|
    \le (d+2)\big(10d(2d-1)M\cL(k-1)\big)^{k},
\end{equation}
i.e. $\cN\in\cS(d,\cL;d+2,10d(2d-1)M)$. Now, Lemma \ref{lem:CompositeDerivative1}
with $\PsiZ(v)=\cL(v-1)$, $g=\ha$, $f=\cN$ and \eqref{bound_N} prove the theorem.
\qed


\begin{remark}
{\rm The arguments from the above proof also imply that \eqref{der-L} holds
for the mixed derivatives of order $k$.
However, Theorem \ref{thm:cutoff2} is sufficient for our purposes in this paper.}
\end{remark}


\begin{remark}
{\rm In Definition \ref{cutoff-second} $\CB_1$ can be replaced by $\CB_\infty$,
but this will lead to some complications in the construction of admissible functions by semi-norms,
as well as bigger constants in Theorem \ref{thm:cutoff2}.}
\end{remark}

For $\cL=\cL_{\ell,\eps}$ the admissible multivariate cutoff function in Theorem \ref{thm:cutoff2} is from
the class $\cS(d,\cL_{\ell,\eps};\gamma_0,\bargamma_0/\eps^2)$,
where the second parameter is of order $\eps^{-2}$ and not of order $\eps^{-1}$ as in the univariate case.
This is due to the method of construction via composition of two functions
from $\cS(1,\cL_{\ell,\eps};\gamma,\bargamma/\eps)$.
The composition necessarily belongs to $\cS(1,\cL_{\ell,\eps};\gamma_0,\bargamma_0/\eps^2)$
unless better estimates for the derivatives are known. A different construction that leads to
a smaller value of the second parameter is given in Subsections \ref{from univariate products}-\ref{small2}.

\subsection{Construction of admissible cutoff functions by univariate products}
\label{from univariate products}

Another natural approach for constructing admissible d-dimensional cutoff functions resembles
the construction of d-dimensional wavelets from univariate father wavelets.


\smallskip

From Lemmas \ref{lem:product_a}-\ref{lem0} we immediately get


\begin{lemma}\label{lem:product_b}
Let $\ha_{1,j}, \ha_{2,j}$, $j=1,\dots,d$, be admissible univariate functions of type $(a)$.
Then
\begin{equation}\label{cutoff-d4b}
    \hB(t)=\prod_{j=1}^d \ha_{1,j}(t_j)-\prod_{j=1}^d \ha_{2,j}(2t_j)
\end{equation}
is an admissible d-dimensional cutoff function of second kind and type $(b)$.
\end{lemma}

In the univariate case all admissible functions of type $(c)$ are among the admissible functions
of type $(b)$ constructed via \eqref{cutoff-d4b}.
Unfortunately, in dimensions $d\ge2$ representation \eqref{cutoff-d4b} does not provide
any admissible d-dimensional function of type $(c)$. In order to get such cutoff functions
we employ two other one-dimensional techniques.


\begin{lemma}\label{lem:product_c}
Let $\hA$ be given by \eqref{cutoff-d3} with $\ha_j$ satisfying $0\le\ha_j(t)\le1$.
We define a cutoff function $\hC$ in two ways, namely,
\begin{equation}\label{cutoff-d5}
    \hC(t):=\begin{cases}
    1-\hA^2(2t),& t\in\CB_\infty,\\
    \hA(t)\sqrt{2-\hA^2(t)},& t\in2\CB_\infty\backslash\CB_\infty,\\
    0, & t\notin2\CB_\infty. \end{cases}
\end{equation}
or
\begin{equation}\label{cutoff-d6}
    \hC(t):=\begin{cases}
    0, &t\in\frac{1}{2}\CB_\infty,\\
    \cos(\frac{\pi}{2}\hA(2t)), &t\in\CB_\infty\backslash\frac{1}{2}\CB_\infty,\\
    \sin(\frac{\pi}{2}\hA(t)), &t\in2\CB_\infty\backslash\CB_\infty,\\
    0,&t\notin2\CB_\infty. \end{cases}
\end{equation}
Then the function $\hC \ge 0$ from \eqref{cutoff-d5} or \eqref{cutoff-d6} is admissible
of second kind and type $(c)$.
\end{lemma}

The proof of this lemma is straightforward.


\begin{remark}
{\rm The cutoff functions constructed in this subsection satisfy a stronger form of
Definition \ref{cutoff-second} with $\CB_1$ replaced by $\CB_\infty$.}
\end{remark}

\subsection{Construction of admissible cutoff functions from univariate products with
``small'' derivatives}\label{small2}

The admissible cutoff functions from univariate products from
\S\ref{from univariate products} allow better estimates on the derivatives than those
in \S\ref{small1}.


\begin{theorem}\label{thm:cutoff2a}
Let $\cL$ satisfy \eqref{def-M-conv} and let $M$ be given by \eqref{def-M}.
Let $\gamma,\bargamma>0$ be such that  the set $\cS(1,\cL;\gamma,\bargamma M)$
contains an admissible univariate cutoff function $\ha$, $0\le\ha(t)\le1$, of type $(a)$ according to
Theorem \ref{thm:cutoff-d=1}.

\smallskip
$(a)$ If $\hA$ is defined by \eqref{cutoff-d3} with $\ha_j=\ha$, then $\hA$ is an admissible
cutoff function of second kind and type $(a)$ and $\hA\in\cS(d,\cL;\gamma,\bargamma M)$.

\smallskip

$(b)$ If $\hB$ is defined by \eqref{cutoff-d4b}
with $\ha_{1,j}=\ha_{1,j}=\ha$, then $\hB$
is an admissible cutoff function of second kind and type $(b)$ and
$\hB\in\cS(d,\cL;2\gamma,2\bargamma M)$.

\smallskip

$(c)$ If $\hC$ is defined by \eqref{cutoff-d6} with $\hA$ from \eqref{cutoff-d3},
where $\ha_j=\ha$, then $\hC$
is admissible of second kind and type $(c)$ and
$\hC\in\cS(d,\cL;1,(\pi\gamma+2)\bargamma M)$.

\end{theorem}

\noindent
{\bf Proof.} It is established in
Lemmas \ref{lem:product_a}, \ref{lem:product_b} and \ref{lem:product_c}
that $\hA, \hB, \hC$ are admissible cutoff functions of the respective type.
The fact that $\hA\in\cS(d,\cL;\gamma,\bargamma M)$ follows
immediately by \eqref{cutoff-d3}, $0\le\ha(t)\le1$ and Definition \ref{Def-S}.
Also $\hB\in\cS(d,\cL;2\gamma,2\bargamma M)$ follows by \eqref{cutoff-d4b}
and Definition \ref{Def-S},
as the constant $2\bargamma$ replaces $\bargamma$
because of the multiplier $2$ in the arguments of the functions in the second product  in \eqref{cutoff-d4b}.

To find bounds on the derivatives of $\hC(t)$ for $1/2\le \|t\|_\infty\le 2$ we fix $1\le j \le d$.
Consider $\hC(t)=(g\circ f)(t_j)$ as function of $t_j\in[1/2,1]$,
where for $\|t\|_\infty\le 1$ we set $f(t_j)=\ha(2t_j)$
and $g(x)=\cos(\lambda x)$ with $\lambda=\frac{\pi}{2}\prod_{m=1,m\ne j}^d\ha(2t_m)$
and for $1<\|t\|_\infty\le 2$ we set $f(t_j)=\ha(t_j)$ and $g(x)=\sin(\lambda x)$
with $\lambda=\frac{\pi}{2}\prod_{m=1,m\ne j}^d\ha(t_m)$.
We apply Lemma \ref{lem:CompositeDerivative1} as
$f\in\cS(1,\cL;\gamma,2\bargamma M)$, $g\in\cS(1,\cL;1,\pi/2)$
and get $g\circ f\in\cS(1,\cL;1,(\pi\gamma+2)\bargamma M)$.
For $t_j\in[1,2]$ we use $f(t_j)=\ha(t_j)$ and $g(x)=\sin(\lambda x)$
with $\lambda=\frac{\pi}{2}\prod_{m=1,m\ne j}^d\ha(t_m)$.
We apply Lemma~\ref{lem:CompositeDerivative1} as
$f\in\cS(1,\cL;\gamma,\bargamma M)$, $g\in\cS(1,\cL;1,\pi/2)$
and get $g\circ f\in\cS(1,\cL;1,(\pi\gamma/2+1)\bargamma M)$.
Consequently, in all cases $\hC\in\cS(d,\cL;1,(\pi\gamma+2)\bargamma M)$.
\qed

\begin{remark}
{\rm In cases (a) and (b) of Theorem \ref{thm:cutoff2a} it sufices to require $\cL$
to satisfy  \eqref{def-M} instead of \eqref{def-M-conv}.}
\end{remark}

\section{Localized tensor product Jacobi polynomial kernels}
\label{d-Jacobi-kernels}
\setcounter{equation}{0}

Denote by $\tP_n^{(\a_j, \b_j)}$ ($1\le j \le d$) the $n$th degree Jacobi polynomial normalized
in $L^2([-1, 1], w_{\a_j, \b_j})$, see \S\ref{Jacobi-kernels}.
Then for multi-indexes $\a=(\a_1, \dots, \a_d)$ and $\b=(\b_1,\dots,\b_d)$
the d-dimensional tensor product Jacobi polynomials are defined by
\begin{equation} \label{product-Jacobi}
\tP_\nu^{(\a, \b)}(x):=\prod_{j=1}^d \tP_{\nu_j}^{(\a_j, \b_j)}(x_j).
\end{equation}
Recall our standing assumption: $\alpha_j, \beta_j \ge -1/2$.
Evidently, $\{\tP_\nu^{(\a, \b)}\}_{\nu\in \NN_0^d}$ is an orthonormal basis for
$L^2([-1, 1]^d, w_{\a, \b})$ with $w_{\a, \b}$ being the product Jacobi weight
defined in \eqref{Jacobi-weight}.

We are interested in kernels of the form
\begin{equation} \label{Jacobi-kernel-d}
\Lambda_n(x, y):= \sum_{\nu\in\NN_0^d} \hA\Big(\frac{\nu}{n}\Big)
\tP_\nu^{(\a, \b)}(x)\tP_\nu^{(\a, \b)}(y),
\quad x, y\in [-1, 1]^d.
\end{equation}

Define
\begin{equation}\label{Jacobi-weight-d}
\WW(n;x):= 
\prod_{j=1}^d w_{\a_j,\b_j}(n; x_j),
\end{equation}
where $w_{\a_j,\b_j}(n; x_j)$ is given in \eqref{Jacobi-weight-n}.
We shall also use the distance on $[-1,1]^d$ defined by
\begin{equation}\label{def-rho-d}
\rho(x, y)=\max_{1\le j\le d}|\arccos x_j-\arccos y_j|.
\end{equation}


\begin{theorem}\label{thm:product-Jacobi-0}
Suppose $\hA \in C^{3k-1}[0, \infty)^d$ for some $k\ge 1$, $\supp \hA \subset [0, 2]^d$,
and for any $t\in [0, 2]^d$ of the form
$t=(t_1, \dots, t_{\ell-1}, 0, t_{\ell+1}, \dots, t_d)$, $1\le \ell\le d$, i.e. $t= \proj_\ell t$,
$\hA$ satisfies
$D^m_\ell\hA(t)=0$ for $m=1, 2, \dots, 3k-1$.
Then the kernels from \eqref{Jacobi-kernel-d} satisfy
\begin{equation}\label{Jacobi-bound2}
|\Lambda_n(x,y)|\le\frac{c n^d }{\sqrt{\WW(n; x)}\sqrt{\WW(n; y)}}\big(1
+ n\rho(x, y)\big)^{-k},
\quad x, y\in [-1, 1]^d.
\end{equation}
Here the constant $c$ depends on $k$, $d$, $\alpha$, $\beta$ and
$\|D^{3k-1}_\ell \hA\|_\infty$, $\ell=1,\dots,d$, 
 but not on $x$, $y$ and $n$.
Consequently, for an admissible cutoff function $\hA$ the above estimate holds
for any $k>0$.
\end{theorem}

\noindent
{\bf Proof.}
Without loss of generality we may assume that $\rho(x, y)=|\arccos x_d-\arccos y_d|$.
We write $\Lambda_n$ from \eqref{Jacobi-kernel-d} as
\begin{multline} \label{Jacobi-kernel-d_2}
\Lambda_n(x, y)
= \sum_{\nu_1=0}^{2n-1}\cdots\sum_{\nu_{d-1}=0}^{2n-1}
\Big[\sum_{\nu_d=0}^{\infty} \hA\Big(\frac{\nu_1}{n},\dots,\frac{\nu_d}{n}\Big)
\tP_{\nu_d}^{(\a_d, \b_d)}(x_d)\tP_{\nu_d}^{(\a_d, \b_d)}(y_d)\Big]\\
\times\prod_{j=1}^{d-1} \tP_{\nu_j}^{(\a_j, \b_j)}(x_j)
\prod_{j=1}^{d-1} \tP_{\nu_j}^{(\a_j, \b_j)}(y_j).
\end{multline}
For any $\nu_1,\dots, \nu_{d-1}$ we estimate the inner sum in \eqref{Jacobi-kernel-d_2}
by using Theorem~\ref{thm:Jacobi-localization-0}.
We get
\begin{multline}\label{Est-Lam-n}
\Big|\sum_{\nu_d=0}^{\infty} \hA\Big(\frac{\nu_1}{n},\dots,\frac{\nu_d}{n}\Big)
\tP_{\nu_d}^{(\a_d, \b_d)}(x_d)\tP_{\nu_d}^{(\a_d, \b_d)}(y_d)\Big|\\
\le \frac{c n}{\sqrt{w_{\a_d,\b_d}(n; x_d)}
\sqrt{w_{\a_d,\b_d}(n; y_d)}}\big(1 + n\rho(x, y)\big)^{-k}.
\end{multline}
For the Jacobi polynomials from the outer products we apply \eqref{est.Pn2},
$\a_j, \b_j\ge -1/2$, and use that $\nu_j<2n$ to obtain
\begin{equation}\label{est-Pnuj}
|\tP_{\nu_j}^{(\a_j,\b_j)}(t)|\le \frac{c}{\sqrt{w_{\a_j,\b_j}(\nu_j; t)}}
\le \frac{c}{\sqrt{w_{\a_j,\b_j}(n; t)}},\quad t=x_j,y_j.
\end{equation}
Combining the above two estimates and the fact that the total number of terms
in the outer sums in (\ref{Jacobi-kernel-d_2}) is $(2n)^{d-1}$ proves the theorem.
$\qed$

We next show that estimate \eqref{Jacobi-bound2} can be improved for
cutoff functions of ``small" derivatives given by
Theorem~\ref{thm:cutoff2} or Theorem~\ref{thm:cutoff2a}.


\begin{theorem}\label{thm:product-Jacobi-1}
Let $\hA$ be an admissible cutoff function
which belongs to $\cS(d,\cL;\gamma,\bargamma M)$
with $\cL$ and $M$ as in \eqref{def-M} and $\gamma,\bargamma>0$
$($see Definition~\ref{Def-S}$)$.
Then the kernels from \eqref{Jacobi-kernel-d} satisfy
\begin{equation} \label{Jacobi-kernel-est}
|\Lambda_n(x, y)|
\le \frac{c n^d }{\sqrt{\WW(n; x)}\sqrt{\WW(n; y)}}
\exp\Big\{-\frac{\tilde{c}n\rho(x, y)}{\cL(n\rho(x, y))}\Big\}
\end{equation}
 for $x, y\in [-1, 1]^d$.
Here $\tilde{c}= c'/\bargamma M$ with $c'>0$ being an absolute constant and
the constant $c>0$ depends on $d$, $M$, $\alpha$, $\beta$, $\gamma$ and $\bargamma$,
 but not on $x$, $y$ and $n$.
\end{theorem}

The proof of Theorem \ref{thm:product-Jacobi-1} is the same as the proof of
Theorem \ref{thm:product-Jacobi-0} with the role
of Theorem \ref{thm:Jacobi-localization-0}
played by Theorem \ref{thm:Jacobi-localization-1} and Remark \ref{rem:Jacobi-localization-1}.

\smallskip

The next theorem shows that the kernels $\Lambda_n(x, y)$ from (\ref{Jacobi-kernel-d})
are ${\rm Lip} 1$ in $x$ and $y$ with respect to the distance $\rho(\cdot,\cdot)$;
it is needed for our further development.


\begin{theorem}\label{thm:Lip}
Under the hypotheses of Theorem~\ref{thm:product-Jacobi-0} with
$
k> 2\max_i\{\a_i+\b_i\}+5
$
for all $x, y, \xi\in [-1, 1]^d$ such that $\rho(x, \xi) \le c_*n^{-1}$, $n\ge 1$, $c_*>0$,
the kernel $\Lambda_n$ from $(\ref{Jacobi-kernel-d})$ satisfies
\begin{equation}\label{Lip}
|\Lambda_n(x,y)-\Lambda_n(\xi,y)|
\le\frac{c n^{d+1}\rho(x, \xi)}{\sqrt{\WW(n; x)}\sqrt{\WW(n; y)}}
\big(1+ n\rho(x, y)\big)^{-\sigma},
\end{equation}
where $\sigma = k - 2\max_i\{\a_i+\b_i\}-5$
and $c>0$ depends only on $k, d, \a, \b$, $c_*$, and 
$\|D^{3k-1}_\ell \hA\|_\infty$, $\ell=1,\dots,d$.
Therefore, for an admissible cutoff function $\hA$ the above estimate holds for any $\sigma>0$.
\end{theorem}

\noindent
{\bf Proof.}
Apparently it suffices to prove estimate \eqref{Lip} for all $\xi\in [-1, 1]^d$ of the form
$\xi=x+\delta e_i$ such that $\rho(x, x+\delta e_i)\le c_*n^{-1}$ and $1\le i\le d$
with $e_i$ being the $i$th coordinate vector.

As in the proof of Theorem~\ref{thm:product-Jacobi-0},
without loss of generality
we may assume that $\rho(x, y)=|\arccos x_d-\arccos y_d|=: \rho_d(x, y)$.
Assuming that $\xi=x+\delta e_i$ is as above, we consider two cases for $i$.


{\em Case 1}: $i=d$. Then we have
\begin{multline*}
\Lambda_n(x, y)-\Lambda_n(x+\delta e_d, y)\\
= \sum_{\nu_1=0}^{2n-1}\cdots\sum_{\nu_{d-1}=0}^{2n-1}
\Big[\sum_{\nu_d=0}^{\infty}\hA\Big(\frac{\nu_1}{n},\dots,\frac{\nu_d}{n}\Big)
\Big(\tP_{\nu_d}^{(\a_d, \b_d)}(x_d)-\tP_{\nu_d}^{(\a_d, \b_d)}(x_d+\delta)\Big)
\tP_{\nu_d}^{(\a_d, \b_d)}(y_d)\Big]\\
\times\prod_{j=1}^{d-1} \tP_{\nu_j}^{(\a_j, \b_j)}(x_j)
\prod_{j=1}^{d-1} \tP_{\nu_j}^{(\a_j, \b_j)}(y_j).
\end{multline*}
Applying Theorem~\ref{thm:Jacobi-Lip-1} to the inner sum we get
\begin{multline*}
\Big|\sum_{\nu_d=0}^{\infty} \hA\Big(\frac{\nu_1}{n},\dots,\frac{\nu_d}{n}\Big)
\Big(\tP_{\nu_d}^{(\a_d, \b_d)}(x_d)-\tP_{\nu_d}^{(\a_d, \b_d)}(x_d+\delta)\Big)
\tP_{\nu_d}^{(\a_d, \b_d)}(y_d)\Big|\\
\le \frac{c n^2\rho(x, x+\delta e_i)}{\sqrt{w_{\a_d,\b_d}(n; x_d)}
\sqrt{w_{\a_d,\b_d}(n; y_d)}}\big(1 + n\rho_d(x, y)\big)^{-\sigma}.
\end{multline*}
For the Jacobi polynomials $\tP_{\nu_j}^{(\a_j, \b_j)}(x_j)$ and
$\tP_{\nu_j}^{(\a_j, \b_j)}(y_j)$
from the outer products we apply estimates \eqref{est-Pnuj}
and combining these with the above we arrive at \eqref{Lip}.

\smallskip


{\em Case 2}: $i\ne d$. Let $x, x+\delta e_i \in[-1, 1]^d$ and
$\rho(x, x+\delta e_i)\le c_*n^{-1}$.
We have
\begin{align}
\Lambda_n(x, y)&-\Lambda_n(x+\delta e_i, y) \notag\\
&= \sum_{\nu_1=0}^{2n-1}\cdots\sum_{\nu_{d-1}=0}^{2n-1}
\Big[\sum_{\nu_d=0}^{\infty}\hA\Big(\frac{\nu_1}{n},\dots,\frac{\nu_d}{n}\Big)
\tP_{\nu_d}^{(\a_d, \b_d)}(x_d)\tP_{\nu_d}^{(\a_d, \b_d)}(y_d)\Big]\label{Lam-n-Lam-n}\\
&\times\prod_{j=1, \, j\ne i}^{d-1} \tP_{\nu_j}^{(\a_j, \b_j)}(x_j)
\Big(\tP_{\nu_i}^{(\a_i, \b_i)}(x_i) - \tP_{\nu_i}^{(\a_i, \b_i)}(x_i+\delta)\Big)
\prod_{j=1}^{d-1} \tP_{\nu_j}^{(\a_j, \b_j)}(y_j).\notag
\end{align}
As is well known that
$
\frac{d}{dt}[P_m^{(\a,\b)}(t)] = \frac{m+\a+\b+1}{2}P_{m-1}^{(\a+1,\b+1)}(t)
$
(see \cite[(4.21.7)]{Sz}).
Combining this with estimate \eqref{est.Pn2} from Lemma~\ref{lem:simple-est}
and
$h_m^{(\a,\b)}\sim h_{m-1}^{(\a,\b)}\sim m^{-1}$ (see \eqref{def-hn})
give
$$
\Big|\frac{d}{dt}\tP_m^{(\a,\b)}(t)\Big|
\le \frac{cm}{\sqrt{w_{\a+1,\b+1}(m-1, t)}}
\le \frac{cm}{\sqrt{w_{\a,\b}(m, t)}\big(\sqrt{1-t^2}+ m^{-1}\big)}.
$$
We use this to obtain for $\theta, \theta' \in [0, \pi]$ with
$|\theta-\theta'|\le c_*m^{-1}$, $m\ge 2$,
\begin{align}\label{Pn-Pn}
|\tP_m^{(\a,\b)}(\cos \theta)-\tP_m^{(\a,\b)}(\cos \theta')|
&\le \frac{cm|\cos \theta - \cos\theta'|}
{\sqrt{w_{\a,\b}(m, \cos \theta)}\big(\sin \theta + m^{-1}\big)}\\
\le \frac{cm\sin|\frac{\theta - \theta'}{2}|\sin|\frac{\theta + \theta'}{2}|}
{\sqrt{w_{\a,\b}(m, \cos \theta)}\big(\sin \theta + m^{-1}\big)}
&\le \frac{cm|\theta - \theta'|}
{\sqrt{w_{\a,\b}(m, \cos \theta)}}.\notag
\end{align}
Note that \eqref{Pn-Pn} is trivial for $m=0,1$.
Therefore,
\begin{equation} \label{P-nu-i-Pnu-i}
|\tP_{\nu_i}^{(\a_i, \b_i)}(x_i) - \tP_{\nu_i}^{(\a_i, \b_i)}(x_i+\delta)|
\le \frac{c\nu_i\rho(x, x+\delta e_i)}{\sqrt{w_{\a_i,\b_i}(\nu_i, x_i)}}
\le \frac{cn\rho(x, \xi)}{\sqrt{w_{\a_i,\b_i}(n, x_i)}}.
\end{equation}
Now, we use \eqref{Est-Lam-n} to estimate the inner sum in \eqref{Lam-n-Lam-n},
\eqref{est-Pnuj} to estimate the Jacobi polynomials $\tP_{\nu_j}^{(\a_j, \b_j)}(x_j)$
($j\ne i$) and $\tP_{\nu_j}^{(\a_j, \b_j)}(y_j)$
from the outer products in \eqref{Lam-n-Lam-n},
and we also use \eqref{P-nu-i-Pnu-i} to obtain again \eqref{Lip}.
Here as well as in Case 1 we took into account that the number of terms in the outer
sums is $(2n)^{d-1}$.
$\qed$

\smallskip

Lower bound estimates for the $L^p$-norms of the kernels $\Lambda_n(x, y)$
in $x$ or $y$ can also be easily derived from the corresponding results
in dimension one.


\begin{proposition}\label{prop:lower-bound}
Let $\hA$ be admissible and $|\hA(t)| \ge c > 0$ for $t \in [1, 1+\delta]^d$, $\delta>0$.
Then for $n\ge 1/\delta$
\begin{equation} \label{est-Lp-norm}
\int_{[-1,1]^d} |\Lambda_n (x,y)|^2 w_{\a,\b}(y) dy \ge c n^d \WW(n;x)^{-d},
\quad x \in [-1,1]^d,
\end{equation}
where $c>0$ depends only on $\delta$, $\a$, $\b$, and $d$.
\end{proposition}

\noindent
{\bf Proof.}
By the definition of $\Lambda_n(x, y)$ in (\ref{Jacobi-kernel-d})
and the orthogonality of the Jacobi polynomials, it follows that
\begin{align*}
& \int_{[-1,1]^d} |\Lambda_n(x,y)|^2 w_{\a,\b}(y) dy   =
    \sum_{\nu \in \NN_0^d} |\hA(\nu/n)|^2  [ \tilde P_\nu^{(\a,\b)}(x)]^2 \\
& \qquad \ge  \sum_{\nu \in [n, n+\delta n]^d} |\hA(\nu/n)|^2  [ \tilde P_\nu^{(\a,\b)}(x)]^2
   \ge  c  \prod_{i=1}^d \sum_{\nu_i = n}^{n+ \lfloor \delta n\rfloor}
         [\tilde P_{\nu_i}^{(\a_i,\b_i)}(x_i)]^2
\end{align*}
and the stated lower bound follows from the respective result in the univariate case,
given in \cite[Proposition 2.4]{KPX1}.
$\qed$

The rapidly decaying polynomial kernels $\Lambda_n(x, y)$ from (\ref{Jacobi-kernel-d})
can be utilized as in the univariate case \cite[Proposition 2.6]{KPX1}
for establishing Nikolski type inequalities:


\begin{proposition}\label{Nikolski}
For $0 < q \le p \le \infty$ and $g \in \Pi_n^d$,
\begin{equation}\label{norm-relation}
\|g\|_p \le cn^{(2d + 2 \sum_{i=1}^d \min\{0, \max\{\a_i, \b_i\}\} )(1/q-1/p)}\|g\|_q,
\end{equation}
furthermore, for any $s\in\RR$,
\begin{equation}\label{norm-relation2}
\|\WW(n;\cdot)^{s} g(\cdot)\|_p
\le cn^{d(1/q-1/p)}\|\WW(n;\cdot)^{s+1/p-1/q} g(\cdot)\|_q.
\end{equation}
\end{proposition}

\section{Additional auxiliary results}\label{background}
\setcounter{equation}{0}

\subsection{The maximal inequality}\label{max}

We let $\cM_t$ $(0<t<\infty)$ be the maximal operator defined by
\begin{equation}\label{def.max-fun}
\cM_tf(x):=\sup_{I\ni x}\left(\frac1{\mu(I)}\int_I|f(y)|^t w_{\a,\b}(y)\, dy\right)^{1/t},
\quad x\in[-1,1]^d,
\end{equation}
where the sup is over all boxes (rectangles) $I\subset[-1,1]^d$ with sides parallel
to the coordinate axces containing $x$.
Here $\mu(E):= \int_E w_{\a,\b}(y)\, dy$.

We denote by $B(\xi, r)$ the ``ball" (box) centered at $\xi\in [-1, 1]^d$ of radius $r>0$
with respect to the distance $\rho(\cdot,\cdot)$ on $[-1, 1]^d$, i.e.
\begin{equation}\label{def-ball}
B(\xi, r):=\{x\in [-1, 1]^d: \rho(x, \xi)<r\}.
\end{equation}
We next show that for $0<\delta \le \pi$
\begin{equation}\label{mu-B}
\mu(B(y, \delta))\sim \delta^d\prod_{i=1}^d \Big(\sqrt{1-y_i^2}+\delta\Big)^{2\gamma_i+1},
\quad
\gamma_i:=
\left\{
\begin{array}{lcl}
\a_i &\mbox{if}& 0\le y_i\le 1,\\
\b_i &\mbox{if}& -1\le y_i<0.
\end{array}
\right.
\end{equation}
Let $y_i=:\cos \phi_i$, $0\le \phi_i\le \pi$, and
$\varphi'_i:=\max\{\phi_i-\delta, 0\}$,
$\varphi''_i:=\min\{\phi_i+\delta, \pi\}$.
Evidently
\begin{align*}
\mu(B(y, \delta))
&= \prod_{i=1}^d \int_{\cos \varphi''_i}^{\cos \varphi'_i} (1-x_i)^{\a_i}(1+x_i)^{\b_i} dx_i\\
& = \prod_{i=1}^d \int_{\varphi'_i}^{\varphi''_i}
(1-\cos \theta_i)^{\a_i}(1+\cos \theta_i)^{\b_i} \sin \theta_i d\theta_i\\
&\sim  \delta^d\prod_{i=1}^d (\sin\phi_i+\delta)^{2\gamma_i+1}
=  \delta^d\prod_{i=1}^d \Big(\sqrt{1-y_i^2}+\delta \Big)^{2\gamma_i+1},
\end{align*}
which confirms (\ref{mu-B}).

By (\ref{mu-B}) it follows that
$\mu(B(y, 2\delta)) \le c\mu(B(y, \delta))$,
i.e. $\mu$ is a doubling measure on $[-1,1]^d$ and,
therefore, the Fefferman-Stein vector-valued maximal inequality is valid
(see \cite{Stein}):
Assuming that $0<p<\infty, 0<q\le\infty$ and $0<t<\min\{p,q\}$,
then for any sequence of functions $\{f_k\}_{k=1}^\infty$ on $[-1,1]^d$,
\begin{equation}\label{max_ineq}
\nnorm{\Bigl(\sum_{k=1}^\infty|\cM_tf_k(\cdot)|^q\Bigr)^{1/q}}_{\Lp}
\le c\nnorm{\Bigl(\sum_{k=1}^\infty| f_k(\cdot)|^q\Bigr)^{1/q}}_{\Lp}.
\end{equation}

We need to estimate $(\cM_t\ONE_{B(y, \delta)})(x)$. Such estimates readily follow
by (\ref{mu-B}) and the respective univaruate result in \cite[Lemma 2.7]{KPX1}.


\begin{lemma}\label{lem:J-maximal}
Let $y \in [-1, 1]^d$ and $0<r\le \pi$, and suppose $\gamma_i$, $i=1, \dots, d$,
are defined as in \eqref{mu-B}.
Then for any $x\in [-1, 1]^d$
\begin{equation}\label{J-max1}
 (\cM_t \ONE_{B(y,r)})(x)
\sim \prod_{j=1}^d \Big(1+\frac{\rho(y_j, x_j)}{r}\Big)^{-1/t}
     \Big(1+\frac{\rho(y_j, x_j)}{r+\rho(y_j, 1)}\Big)^{-(2\gamma_j+1)/t}
\end{equation}
and hence
\begin{equation}\label{J-max3}
(\cM_t \ONE_{B(y, r)})(x)
\ge c \prod_{j=1}^d \Big(1+\frac{\rho(y_j, x_j)}{r}\Big)^{-(2\gamma_j+2)/t}
\ge c  \Big(1+\frac{\rho(y, x)}{r}\Big)^{-(2|\gamma|+2d)/t}.
\end{equation}
Here $\rho(y_j, x_j):=|\arccos y_j-\arccos x_j|$
and $\rho(y, x)$ is defined in $(\ref{def-rho-d})$.
\end{lemma}

We also want to record the following useful inequality which follows easily from the case $d=1$,
proved in \cite[(2.22)]{KPX1}:
\begin{equation}\label{omega<omega}
\WW(n; x)\le  c\WW(n; y)(1+n \rho(x, y))^{d+ 2\sum_{i=1}^d \max\{\a_i, \b_i\}},
\end{equation}
for  $x, y\in [-1, 1]^d$ and $n\ge 1$,
where $\WW(n; x)$ is from \eqref{Jacobi-weight-d}.

\subsection{Distributions on \boldmath $[-1, 1]^d$}\label{distributions}

Here we introduce and give some basic facts about distributions on $[-1, 1]^d$.
We shall use as test functions the set $\cD:=C^{\infty}[-1,1]^d$,
where the topology is induced by the semi-norms
\begin{equation}\label{semi-norms}
|\phi|_\mu:= \| D^\mu \phi(t)\|_\infty
\quad \hbox{for all multi-indices $\mu$.}
\end{equation}
Observe that the tensor product Jacobi polynomials $\{\tilde P_\nu^{(\a,\b)}\}$ belong to $\cD$
and more importantly the test functions $\phi\in\cD$
can be completely characterized by the coefficients of their Jacobi expansions.
Denote
\begin{equation}\label{D-norms}
\cN_k(\phi):=\sup_{\nu \in \NN_0^d}\, (|\nu|+1)^k |\langle \phi, \tilde P_\nu^{(\a,\b)}\rangle|,
\end{equation}
where $\langle f, g \rangle :=\int_{[-1,1]^d} f(x)\overline{g(x)} w_{\a,\b}(x)dx$.


\begin{lemma}\label{lem:char-D}
$(i)$
$\phi\in\cD$ if and only if
$|\langle \phi, \tilde P_\nu^{(\a,\b)}\rangle| =\cO((|\nu|+1)^{-k})$ for all $k$.

$(ii)$
For every $\phi\in \cD$ we have
$
\phi=\sum_{\nu \in \NN_0^d} \langle \phi, \tilde P_\nu^{(\a,\b)}\rangle\tilde P_\nu^{(\a,\b)},
$
where the convergence is in the topology of $\cD$.

$(iii)$
The topology in $\cD$ can be equivalently defined by the norms
$\cN_k(\cdot)$, $k\ge 0$.
\end{lemma}
The proofs of this lemma is easy and similar to the proof of Lemma~2.8 in \cite{KPX1}.

\smallskip

The space $\cD'$ of distributions on $[-1,1]^d$ is defined as the set
of all continuous linear functionals on $\cD$.
The pairing of $f\in \cD'$ and $\phi\in\cD$ will usually be denoted by
$\langle f, \phi \rangle:= f(\overline{\phi})$.
As~will be shown it is in a sense consistent with the inner product
$
\langle f, g \rangle  
$
in $L^2(w_{\a,\b})$.
We shall need the representation of distributions from $\cD'$
in terms of Jacobi polynomials.


\begin{lemma}\label{lem:dec-D1}
$(i)$
A linear functional $f$ on $\cD$ belongs to $\cD'$
if and only if there exists $k\ge 0$ such that
\begin{equation}\label{D1}
|f(\phi)|=|\langle f, \overline\phi\rangle|\le c_k\cN_k(\phi)
\quad \mbox{for all } \; \phi \in \cD,
\end{equation}

$(ii)$
For any $f\in\cD'$ there exist constants $c>0$ and $k\ge 0$ such that
\begin{equation}\label{D2}
|f(\tP_\nu^{(\a,\b)})|=|\langle f, \tP_\nu^{(\a,\b)} \rangle| \le c_k(|\nu|+1)^k
\quad \mbox{for all } \;  \nu \in \NN_0^d, \quad \mbox{and}
\end{equation}
\begin{equation}\label{D3}
f(\phi)=\lim_{n\to\infty} \langle S_n, \overline{\phi}\rangle
=\sum_{\nu \in \NN_0^d}
\langle f, \tilde P_\nu^{(\a,\b)}\rangle \langle \phi, \tilde P_\nu^{(\a,\b)}\rangle
\quad \mbox{for} \quad \phi\in \cD,
\end{equation}
where
$
S_n:=\sum_{|\nu|\le n} \langle f, \tilde P_\nu^{(\a,\b)}\rangle\tilde P_\nu^{(\a,\b)}
$
and the series converges absolutely.

$(iii)$
For any sequence $\{c_\nu\}_{\nu\in\NN_0^d}$ satisfying
$|c_\nu| \le A(|\nu|+1)^\ell$ for $\nu\in\NN_0^d$ and some constants $A$ and $\ell$,
the sequence
$$
s_n:= \sum_{|\nu|\le n} c_\nu P_\nu^{(\a,\b)}
$$
converges in $\cD'$ as $n \to \infty$ to some distribution $F\in\cD'$ such that
$\langle F, \tP_\nu^{(\a,\b)} \rangle = c_\nu$ for $\nu\in\NN_0^d$.

\end{lemma}

\noindent
{\bf Proof.}
Part (i) of the lemma follows by the definition of $\cD'$ and Lemma~\ref{lem:char-D}
as in the classical case.

Estimate \eqref{D2} is immediate from \eqref{D1} and \eqref{D-norms}.
Further, we have for $\phi\in\cD$
$$
\lim_{n\to\infty} \langle S_n, \overline{\phi} \rangle
= \lim_{n\to\infty} f\Big(\sum_{|\nu|\le n} \langle \phi, \tilde P_\nu^{(\a,\b)}\rangle P_\nu^{(\a,\b)}\Big)
=f(\phi),
$$
which confirms \eqref{D3}. Here we used Lemma~\ref{lem:char-D}, (ii).

To prove part (iii), we observe that
$
\langle s_n, \overline{\phi} \rangle = \sum_{|\nu|\le n} c_\nu\langle \phi, \tP_\nu^{(\a, \b)} \rangle
$
for $\phi\in\cD$ and using the assumption and Lemma~\ref{lem:char-D} we get
$
|c_\nu| |\langle \phi, \tP_\nu^{(\a, \b)} \rangle| \le c(|\nu|+1)^{\ell-k}
$
for an arbitrary $k\ge 0$.
Therefore, the series $\sum_{\nu\in\NN_0^d} c_\nu\langle \phi, \tP_\nu^{(\a, \b)} \rangle$
converges absolutely and hence
\begin{equation}\label{rep-F-phi}
F(\phi):= \lim_{n\to\infty}\langle s_n, \overline{\phi} \rangle
= \sum_{\nu\in\NN_0^d} c_\nu\langle \phi, \tP_\nu^{(\a, \b)} \rangle,
\quad \phi\in\cD,
\end{equation}
is a well defined linear functional.
We claim that $F$ is bounded. Indeed, for $\phi\in\cD$
\begin{align*}
|F(\phi)|
&\le \sum_{\nu\in\NN_0^d} |c_\nu||\langle \phi, \tP_\nu^{(\a, \b)} \rangle|
\le A\sum_{\nu\in\NN_0^d}(|\nu|+1)^\ell|\langle \phi, \tP_\nu^{(\a, \b)} \rangle|\\
& \le A\cN_{\ell+d+1}(\phi)\sum_{\nu\in\NN_0^d}(|\nu|+1)^{-d-1}
\le c\cN_{\ell+d+1}(\phi),
\end{align*}
which shows that $F\in\cD'$.

Finally,
$F(\tP_\nu^{(\a,\b)})
= \lim_{n\to \infty} \langle s_n, \overline{\phi} \rangle = c_\nu$
is immediate by \eqref{rep-F-phi}.
$\qed$

\smallskip

To simplify our notation, we introduce the following ``convolution": For functions
$\Phi: [-1, 1]^d\times[-1, 1]^d \to \bC$ and $f: [-1, 1]^d \to \bC$, we define
\begin{equation}\label{convolution}
         \Phi*f(x) := \int_{[-1,1]^d} \Phi(x, y)f(y) w_{\a,\b}(y)\,dy
\end{equation}
and extend it to $\cD'$ by duality, i.e.
assuming that $f \in \cD'$ and $\Phi: [-1, 1]^d\times [-1, 1]^d \to\bC$ is such that
$\Phi(x, y)$ belongs to $\cD$ as a function of $y$,
we define $\Phi*f$ by
\begin{equation}\label{convolution1}
\Phi*f (x) := \langle f, \overline{\Phi(x, \cdot)} \rangle.
\end{equation}
Here on the right $f$ acts on $\overline{\Phi(x, y)}$ as a function of $y$.

\subsection{\boldmath $L^p$-multipliers}

We shall need $L^p$-multipliers  for tensor product Jacobi poly-\\nomial expansions.
Since we cannot find any such multipliers in the literature
we next derive simple but non-optimal multipliers satisfying the First Boundary Condition (\S\ref{remedy})
of a certain order.


\begin{theorem}\label{thm:multipliers}
Let $m\in C^r[0, \infty)^d$ for $r$ sufficiently large
$($$r>6\max_i\{\a_i+\b_i\}+6\sum_i\max\{\a_i,\b_i\}+ 6d+20$ will do$)$
and suppose $m$ satisfies the following condition: For any $t\in [0, \infty)^d$ of the form
$t=(t_1, \dots, t_{\ell-1}, 0, t_{\ell+1}, \dots, t_d)$, $1\le \ell\le d$,
we have
$D^s_\ell m(t)=0 $ for $s=1, 2, \dots, r$.
Also, assume
\begin{equation}\label{derivatives}
|D^\tau m(t)| \le c(1+\|t\|_\infty)^{-|\tau|}
\quad \mbox{ for $t\in [0, \infty)^d$ and $|\tau| \le r$, }
\end{equation}
with $c>0$ independent of $t$.
Then the operator
$T_m f := \sum_{\nu\in\NN_0^d} m(\nu) \langle f, \tilde P_\nu^{(\a,\b)}\rangle \tilde P_\nu^{(\a,\b)}$
is bounded on $L^p(w_{\a,\b})$ for $1<p<\infty$.
\end{theorem}

\noindent
{\bf Proof.}
We shall utilize a standard decomposition of unity argument.
Let $\hC$ be an admissible cutoff function of
type (c). Then $\hB=|\hC|^2\ge 0$ is admissible of
type (b) and $ \sum_{j=0}^\infty \hB(2^{-j}t)=1 $ for
$t\in[0, \infty)^d\setminus [0, 1)^d$.
We define $\Phi_0(x,y):=m(0)\tilde P_0^{(\a,\b)}(x)\tilde P_0^{(\a,\b)}(y)$ and
$$
\Phi_j(x, y):= \sum_{\nu\in\NN_0^d} \hB\Big(\frac{\nu}{2^{j-1}}\Big)m(\nu)\tPnu(x)\tPnu(y),
\quad j\ge 1.
$$
Consider the kernels
$K_N:=\sum_{j=0}^N \Phi_j$.
We shall prove that
\begin{equation}\label{est-KN}
\|K_N*f\|_p \le c\|f\|_p
\quad \mbox{for $f\in L^p(w_{\a,\b})$}
\end{equation}
with $c>0$ a constant independent of $f$ and $N$.
As a consequence of this, it is easy to show that
for any  $f\in L^p(w_{\a,\b})$ one has
$T_m f = \lim_{N\to\infty}K_Nf$ in  $L^p(w_{\a,\b})$
and $\|T_mf\|_p \le c\|f\|_p$
as claimed.

To prove (\ref{est-KN}) we shall employ the theory of generalized Calde\'{o}n-Zygmund operators.
Note first that by Parseval's identity
\begin{equation}\label{L2-est}
\|K_N*f\|_2 \le c\|f\|_2
\quad \mbox{for $\;f\in L^2(w_{\a,\b}).$}
\end{equation}
Following Stein \cite{Stein}, p. 29, denote
\begin{equation}\label{def:V}
V(x, y):= \inf \{\mu(B(y, \delta)): x\in B(y, \delta)\}= \mu(B(y, \rho(x, y))),
\end{equation}
where the last equality follows from the definition of $\rho(\cdot, \cdot)$. We shall show that
\begin{equation}\label{est-KN-KN}
|K_N(x, y)-K_N(x, \bar y)| \le c\frac{\rho(y, \bar y)}{\rho(x,\bar y)}[V(x, \bar y)]^{-1}
\end{equation}
whenever
$\rho(x,\bar y) \ge 2\rho(y,\bar y)$.
Then  (\ref{est-KN}) will follow for $1<p\le 2$ by (\ref{L2-est}) and (\ref{est-KN-KN})
using the proposition on pp. 29-30 and Theorem 3 on p. 19 in \cite{Stein}.
After that a standard duality argument leads to estimate (\ref{est-KN})
in the case $2<p<\infty$.

We now turn to the proof of (\ref{est-KN-KN}).
Fix $x, y, \y \in [-1, 1]^d$, $x\ne \y$,
and define
$$
\gamma_i:=
\left\{
\begin{array}{lcl}
\a_i &\mbox{if}& 0\le \y_i\le 1,\\
\b_i &\mbox{if}& -1\le \y_i<0.
\end{array}
\right.
$$
By \eqref{def:V} and \eqref{mu-B} it follows that
\begin{equation}\label{est-V}
V(x, \y)
\sim  \rho(x, \y)^d\prod_{i=1}^d \Big(\sqrt{1-\y_i^2}+ \rho(x, \y)\Big)^{2\gamma_i+1}.
\end{equation}

Let $\hA_j(t):= \hB(t)m(2^{j-1}t)$.
We have $\supp \hB \subset [0, 2]^d\setminus \frac{1}{2}\cB_1$ and by \eqref{derivatives}
it readily follows that
$|D^\tau[m(2^{j-1}t)]| \le c$ for $t\in [0, 2]^d\setminus \frac{1}{2}\cB_1$
and $|\tau|\le r$,
where the constant $c>0$ is independent of $j$.
Therefore,
$\|D^\tau \hA_j\|_\infty = \|D^\tau[\hB(\cdot)m(2^{j-1}\cdot)]\|_\infty \le c$
for $|\tau|\le r$ with $c>0$ independent of $j$.
Now, it is evident that $\hA_j$ satisfies the assumptions of Theorem~\ref{thm:Lip}
for some
$k>2\max_i\{\a_i+\b_i\}+2\sum_i\max\{\a_i,\b_i\}+2d+6$
and hence, using also \eqref{omega<omega}, we get
\begin{equation}\label{Phi-Phi}
|\Phi_j(x, y) - \Phi_j(x, \bar y)|
\le \frac{c2^{j(d+1)}\rho(y, \bar y)}{W_{\a,\b}(2^j, \bar y)(1+2^j\rho(x, \bar y))^\sigma}
\end{equation}
if $\rho(y, \y)\le 2^{-j}$, where $\sigma=k-2\max_i\{\a_i+\b_i\}-5$.

If $\rho(y, \y) > 2^{-j}$ and $\rho(x, \y) \ge 2\rho(y,\y)$ (hence $\rho(x,y)\ge\rho(y,\y)$
and $\rho(x,\y)\le2\rho(x,y)\le3\rho(x,\y)$), then
estimate \eqref{Phi-Phi} follows by Theorem~\ref{thm:product-Jacobi-0} applied separately to
$\Phi_j(x, y)$ and $\Phi_j(x, \bar y)$ and using (\ref{omega<omega}).
Therefore, \eqref{Phi-Phi} holds whenever $\rho(x, \y) \ge 2\rho(y,\y)$.

Let $2^{-j_1-1} \le \rho(x,\y) < 2^{-j_1}$.
Then using $\Phi_0(x, y) = \Phi_0(x, \bar y)$ we write
\begin{align*}
|K_N(x, y) - K_N(x, \bar y)|
&\le \sum_{j=1}^{j_1} |\Phi_j(x, y) - \Phi_j(x, \bar y)|
+ \sum_{j=j_1+1}^N |\Phi_j(x, y) - \Phi_j(x, \bar y)|\\
& =: F_1+F_2.
\end{align*}
%
%
%
For $F_1$ we have using (\ref{Phi-Phi}) and \eqref{est-V}
\begin{align*}
F_1
&\le  \frac{c\rho(y, \bar y)}{\prod_{i=1}^d \big( \sqrt {1-\y_i^2}+2^{-j_1}\big)^{2\gamma_i+1}}
\sum_{j=1}^{j_1} 2^{j(d+1)}\\
&\le  \frac{c\rho(y, \y)2^{j_1(d+1)}}{\prod_{i=1}^d \big( \sqrt {1-\y_i^2}+2^{-j_1}\big)^{2\gamma_i+1}} \\
&\le  \frac{c\rho(y, \y)}{\rho(x, \y)^{d+1}\prod_{i=1}^d \big( \sqrt {1-\y_i^2}+\rho(x, \y)\big)^{2\gamma_i+1}}
\le c\frac{\rho(y, \y)}{\rho(x, \y)}[V(x, \y)]^{-1}.
\end{align*}
%
%
%
To estimate $F_2$ we first observe that
$\big( \sqrt {1-\y_i^2}+2^{-j}\big)(1+2^{j-j_1}) \ge  \sqrt {1-\y_i^2}+2^{-j_1}$.
Then, using again (\ref{Phi-Phi}) and \eqref{est-V}, we get
\begin{align*}
F_2
&\le  c\rho(y, \bar y)
\sum_{j=j_1+1}^N \frac{2^{j(d+1)}}
{\prod_{i=1}^d \big( \sqrt {1-\y_i^2}+2^{-j}\big)^{2\gamma_i+1}  \big(1+2^{j-j_1}\big)^\sigma}\\
&\le  c\rho(y, \bar y)
\sum_{j=j_1+1}^\infty \frac{2^{j(d+1)}}
{\prod_{i=1}^d \big( \sqrt {1-\y_i^2}+2^{-j_1}\big)^{2\gamma_i+1}
\big(1+2^{j-j_1}\big)^{\sigma-2\sum_i\gamma_i-d}}\\
&\le  \frac{c\rho(y, \bar y) 2^{j_1(d+1)}}
{\prod_{i=1}^d \big( \sqrt {1-\y_i^2}+2^{-j_1}\big)^{2\gamma_i+1} }
\sum_{j=j_1+1}^\infty 2^{-(j-j_1)(\sigma-2\sum_i\gamma_i-2d-1)}\\
&\le  \frac{c\rho(y, \bar y)}{\rho(x, \y)^{d+1}
\prod_{i=1}^d \big( \sqrt {1-\y_i^2}+\rho(x,\y)\big)^{2\gamma_i+1}}
\le c\frac{\rho(y, \bar y)}{\rho(x, \y)}[V(x, \y)]^{-1},
\end{align*}
where we used that
$\sigma> 2\sum_i\gamma_i+2d+1$.
The above estimates of $F_1$ and $F_2$ yield (\ref{est-KN-KN}).
This completes the proof of the proposition.
$\qed$


\section{Construction of building blocks (Needlets)}\label{Needlets}
\setcounter{equation}{0}

The construction of frames (needlets) on $[-1,1]^d$ has two basic components:
(i) a~Calder\'{o}n type decomposition formula and
(ii)  a cubature formula.

\subsection{Cubature formula and subdivision of \boldmath $[-1,1]^d$}
\label{sec:cubature}

For the construction of needlets we shall employ the Gaussian quadrature formula
on $[-1, 1]$ with weight
$w_{\a,\b}(t):=(1-t)^\a(1-t)^\b$.
Given $j\ge 0$,
denote by $\xi^m=:\cos \theta_{m}$, $m=1, 2,\dots, 2^{j+1}$,
the zeros of the Jacobi polynomial $P_{2^{j+1}}^{(\a,\b)}$
ordered so that
$
0< \theta_{1} < \dots < \theta_{2^{j+1}} < \pi
$
and set $$\cX_j^{\a,\b}:= \{\xi^m: 1\le m \le 2^{j+1}\}.$$
It is well known that uniformly (see \cite{FW})
\begin{equation}\label{distr-zeros}
\theta_1 \sim 2^{-j},
\;\;
\pi-\theta_{2^{j+1}}  \sim 2^{-j},
\;\;
\theta_{m+1} - \theta_m \sim 2^{-j},
\;\;
\mbox{and hence} \quad \theta_m \sim m2^{-j}.
\end{equation}
As is well known \cite{Sz}
the zeros of the Jacobi polynomial
$P_{2^{j+1}}^{(\a,\b)}$ serve as knots of the Gaussian quadrature
\begin{equation}\label{quadrature-Gauss}
\int_{[-1,1]} f(t)w_{\a,\b}(t)dt \sim \sum_{\xi \in \cX_j^{\a,\b}} \cc_\xi f(\xi),
\end{equation}
which is exact for all algebraic polynomials that are of degree $2^{j+2}-1$.
Furthermore, the coefficients $\cc_\xi$ are all positive and satisfy (see e.g. \cite{N})
\begin{equation}\label{quadrat-coeff}
c_\xi\sim 2^{-j}w_{\a,\b}(\xi)(1-\xi^2)^{1/2}.
\end{equation}

\noindent
{\bf Tiling of \boldmath $[-1, 1]$.}
With $\{\xi^m\}$ as above we write
$$
I_{\xi^m}:= [(\xi^{m+1}+\xi^m)/2, (\xi^{m-1}+\xi^m)/2],
\quad m=2, 3, \dots, 2^{j+1}-1,
$$
and
$$
I_{\xi^1}:= [(\xi^2+\xi^1)/2, 1],
\quad  I_{\xi^{2^{j+1}}}:= [-1, (\xi^{2^{j+1}}+\xi^{2^{j+1}-1})/2].
$$
We define
$$
\cI_j^{\a,\b}:= \{I_{\xi^m}: 1\le m \le 2^{j+1}\}.
$$


For multi-indices $\a=(\a_1, \dots, \a_d)$, $\b=(\b_1, \dots, \b_d)$
and $j\ge 0$, $1\le i\le d$, we denote by
$\cX_j^{\a_i,\b_i}$
the zeroes of the Jacobi polynomial $P_{2^{j+1}}^{(\a_i,\b_i)}$
and write
\begin{equation}\label{def-Xj}
\cX_j:= \cX_j^{\a_1,\b_1}\times\cdots\times \cX_j^{\a_d,\b_d}.
\end{equation}
Now, for $\xi=(\xi_1, \dots, \xi_d)\in \cX_j$ we set
$c_\xi:=c_{\xi_1}\cdots c_{\xi_d}$, where $c_{\xi_i}$ is the corresponding coefficient
of the Gaussian quatrature (\ref{quadrature-Gauss}) with $\a=\a_i$ and $\b=\b_i$.
Evidently, the cubature formula
\begin{equation}\label{quadrature1}
\int_{[-1,1]^d} f(x)w_{\a,\b}(x)dx \sim \sum_{\xi \in \cX_j} \cc_\xi f(\xi)
\end{equation}
is exact for all polynomials in d-variables of degree $2^{j+2}-1$ in each variable
and by (\ref{quadrat-coeff}) the coefficients $\{\cc_\xi\}$ are positive and satisfy
\begin{equation}\label{coeff}
\cc_\xi \sim 2^{- d j} \WW(2^j; \xi),
\end{equation}
where $\WW(2^j; \xi)$ is defined in (\ref{Jacobi-weight-d}).

\medskip

\noindent
{\bf Tiling of \boldmath $[-1, 1]^d$.}
For $\xi=(\xi_1,\dots,\xi_d)\in \cX_j$ , we write
\begin{equation}\label{def-I-xi}
I_\xi : = I_{\xi_1} \times \cdots\times I_{\xi_d},
\quad I_{\xi_i}\in\cI_j^{\a_i,\b_i}.
\end{equation}
Evidently, $[-1, 1]^d=\cup_{\xi\in\cX_j} I_\xi$ and
the interiors of the tiles $\{I_\xi\}_{\xi\in\cX_j}$ do not overlap.

With $B(y,r)$ defined in (\ref{def-ball}) it easily follows from
the univariate case that there exist constants $c_1, c_2 >0$ such that
\begin{equation}\label{I-B}
B(\xi, c_12^{-j}) \subset I_\xi \subset B(\xi, c_22^{-j}),
\quad \xi\in \cX_j.
\end{equation}
By (\ref{mu-B}) it follows that
\begin{equation}\label{muI}
\mu(I_\xi):= \int_{I_\xi}w_{\a,\b}(x)\, dx \sim 2^{-j}\WW(2^j; \xi)
     \sim \cc_\xi,  \quad \xi\in \cX_j,\quad j\ge 0.
\end{equation}

The next lemma is of an independent interest and is instrumental in the subsequent development.


\begin{lemma}\label{lem:quasi-Lip}
Let $P\in \Pi_{2^j}^d$, $j \ge 0$, and $\xi\in \cX_j$.
Suppose $x',x''\in [-1,1]^d$ and
$\rho(x',\xi)\le c_{\star}2^{-j}$,
$\rho(x'',\xi)\le c_{\star}2^{-j}$.
Then for any $\sigma>0$
$$
|P(x')-P(x'')|\le c_\sigma 2^j \rho(x', x'')
\sum_{\eta\in\cX_j}\frac{|P(\eta)|}{(1+2^j \rho(\xi, \eta))^\sigma},
$$
where $c_\sigma>0$ depends only on $\sigma$, $\a$, $\b$, $d$, and $c_\star$.
\end{lemma}

The proof of this lemma is merely a repetition of the proof of the univariate
result in \cite[Lemma 9.2]{KPX1} and will be omitted.

\subsection{Needlets on \boldmath $[-1,1]^d$}\label{sec:Needlets}

The construction of needlet systems is now standard and follows a well established scheme.
We begin with two cutoff functions  $\hA$, $\hB$ of type (b) which satisfy
(see Lemma~\ref{lem:dual-cutoff}):
\begin{equation}\label{part-unity}
  \sum_{j=0}^\infty \overline{\hA(2^{-j}t )}\, \hB(2^{-j}t)=1,
\quad t\in [0,\infty)^d \setminus \mathcal{B}_\infty.
\end{equation}
We define
$\Phi_0(x, y)= \Psi_0(x, y):= \tilde P_0(x)\tilde P_0(y)$,
\begin{align}
\Phi_j(x, y) &:= \sum_{\nu \in \NN_0^d}
\hA\Big(\frac{\nu}{2^{j-1}}\Big)\tilde P_\nu(x)\tilde P_\nu(y),
\quad j\ge 1, \; \mbox{and} \label{def.Phi-j}\\
\Psi_j(x, y) &:= \sum_{\nu \in \NN_0^d}
   \hB \Big(\frac{\nu}{2^{j-1}}\Big)\tilde P_\nu(x)\tilde P_\nu(y),
\quad j\ge 1. \label{def-Psi-j}
\end{align}
Let $\cX_j$ be the set of knots of cubature formula \eqref{quadrature1},
defined in (\ref{def-Xj}), and let $\{c_\xi\}$ be its coefficients.
We define the $j$th level {\em needlets} by
\begin{equation}\label{def-needlets1}
\ph_\xi(x) := \cc_\xi^{1/2}\Phi_j(x, \xi)
\quad\mbox{and}\quad
\psi_\xi(x) := \cc_\xi^{1/2}\Psi_j(x, \xi),
\qquad \xi \in \cX_j.
\end{equation}
We write $\cX := \cup_{j = 0}^\infty \cX_j$,
where equal points from different levels $\cX_j$ are considered
as distinct elements of $\cX$, so that $\cX$ can be used as an index set.
We define the {\em analysis} and {\em synthesis} needlet systems
$\Phi$ and $\Psi$ by
\begin{equation}\label{def-needlets2}
\Phi:=\{\ph_\xi\}_{\xi\in\cX}, \quad \Psi:=\{\psi_\xi\}_{\xi\in\cX}.
\end{equation}

Theorem~\ref{thm:product-Jacobi-0} and (\ref{omega<omega})
imply that the needlets decay rapidly, namely,
\begin{equation}\label{local-needlets3}
  |\ph_\xi(x)|, |\psi_\xi(x)|
\le \frac{c_\sigma2^{jd/2}}{\sqrt{\WW(2^j; \xi)}}\big(1+2^{j}\rho(\xi, x)\big)^{-\sigma},
\quad x\in [-1, 1]^d, \; \forall \sigma.
\end{equation}

We next give estimates on the norms of the needlets, which can be proved
exactly as in the case $d=1$, upon using (\ref{local-needlets3})
and the lower bound estimate from Proposition~\ref{prop:lower-bound}:
For $0<p\le\infty$,
\begin{equation}\label{norm-Needlets}
\|\ph_\xi\|_\Lp \sim \|\psi_\xi\|_\Lp \sim \|\tONE_{I_\xi}\|_\Lp
\sim \Big(\frac{2^{d j}}{\WW(2^j; \xi)}\Big)^{1/2-1/p},
\quad \xi\in\cX_j.
\end{equation}
Here $\tONE_{I_\xi}:=\mu(I_\xi)^{-1/2}\ONE_{I_\xi}$
with $\ONE_E$ being the characteristic function of the set $E$.
Moreover, there exist constants $c^*, c^\diamond >0$ such that
\begin{equation}\label{norm-Needlets2}
\|\ph_\xi\|_{L^\infty(B(\xi, c^*2^{-j}))},\;
\|\psi_\xi\|_{L^\infty(B(\xi, c^*2^{-j}))}
\ge c^\diamond \Big(\frac{2^{d j}}{\WW(2^j; \xi)}\Big)^{1/2}.
\end{equation}


The needlet decomposition of $\cD'$ and $L^p$ follows as in the univariate case
(see \cite[Proposition 3.1]{KPX1}) by the definition of needlets
and their superb localization.

\begin{proposition}\label{prop:needlet-rep}
$(i)$ For $f \in \cD'$, we have
\begin{align}
f &= \sum_{j=0}^\infty \Psi_j*\overline{\Phi}_j*f
\quad\mbox{in}\; \cD',\quad \mbox{and}\label{Needle-rep}\\
f &= \sum_{\xi \in \cX} \langle f, \ph_\xi\rangle \psi_\xi
\quad\mbox{in}\; \cD'. \label{needlet-rep1}
\end{align}

$(ii)$ If $f \in \Lpp$, $1\le p \le \infty$, then
$(\ref{Needle-rep})-(\ref{needlet-rep1})$ hold in $\Lpp$.
Moreover, if $1 < p < \infty$, then the convergence in
$(\ref{Needle-rep})-(\ref{needlet-rep1})$ is unconditional.
\end{proposition}


\begin{remark}\label{rem:needlets}
{\rm
(i) Pick $\hA\ge 0$ a cutoff function of type (c) (see Definition~\ref{cutoff-first}).
Then we can choose $\hB=\hA$ in the constuction of needlets in (\ref{part-unity})-(\ref{def-needlets1})
and obtain $\ph_\xi=\psi_\xi$.
Consequently, (\ref{needlet-rep1}) becomes
$f = \sum_{\xi \in \cX} \langle f, \psi_\xi\rangle \psi_\xi$
and it is easy to prove that (see e.g. \cite{KPX1})
$
\|f\|_2 = \Big(\sum_{\xi\in\cX} |\langle f, \psi_\xi\rangle|^2\Big)^{1/2}
$
for $f\in L^2(w_{\a,\b})$, which shows that $\Psi$ is a tight frame for $L^2(w_{\a,\b})$.

(ii)
If $\hA\ge 0$ is an admissible cutoff function of second kind and type (c)
(see Definition~\ref{cutoff-second})
which belongs to $\cS(d,\cL;\gamma,\bargamma M)$,
then Theorem~\ref{thm:product-Jacobi-1} implies sup-exponential localization of the needlets,
namely,
\begin{equation} \label{localiz-psi}
|\psi_\xi(x)|
\le \frac{c 2^{jd/2} }{\sqrt{\WW(2^j; \xi)}}
\exp\Big\{-\frac{\tilde{c}2^j\rho(\xi, x)}{\cL(2^j\rho(\xi, x))}\Big\},
\quad x\in [-1, 1]^d.
\end{equation}
}
\end{remark}

\section{Weighted Triebel-Lizorkin  spaces on $[-1, 1]^d$}\label{Tri-Liz}
\setcounter{equation}{0}

We next utilize the general idea of using spectral or orthogonal decompositions
(see e.g. \cite{Pee, T1}) to introduce weighted Triebel-Lizorkin spaces on $[-1, 1]^d$.
The~theory of these spaces is entirely parallel to their theory in the univariate case,
developed in \cite{KPX1}. Therefore, we shall only state the main results,
provide the important ingredients and refer the reader to \cite{KPX1} for the proofs.

Given an admissible cutoff function $\hA$ of type (b) (see Definition~\ref{cutoff-first})
satisfying the dyadic covering condition
(\ref{star}) we define a sequence of kernels $\{\Phi_j\}$ by
$\Phi_0(x, y) := \tilde P_0(x)\tilde P_0(y)$ and
\begin{equation}\label{def-Phi-j}
\Phi_j(x, y) := \sum_{\nu \in \NN_0^d}
\hA \Big(\frac{\nu}{2^{j-1}}\Big)\tilde P_\nu(x)\tilde P_\nu(y),
\quad j\ge 1.
\end{equation}


\begin{definition}
For $s, \r \in \RR$, $0<p<\infty$, and $0<q\le\infty$
the weighted Triebel-Lizorkin space $\Fsrpq:=\Fsrpq(w_{\a,\b})$ is defined as
the set of all $f\in \cD'$ such that
\begin{equation}\label{Tri-Liz-norm}
\|f\|_{\Fsrpq}:=\Big\|\Big(\sum_{j=0}^{\infty}
\Big[2^{sj}\WW(2^j;\cdot)^{-\r/d}|\Phi_j*f(\cdot)|\Big]^q\Big)^{1/q}\Big\|_{\Lp} <\infty
\end{equation}
with the usual modification when $q=\infty$.
\end{definition}

Note that the above definition is independent of the choice
of $\hA$ as long as $\hA$ is an admissible function of type (b),
satisfying (\ref{star}) (see Theorem~\ref{thm:Fnorm-equivalence} below).

Also, $\Fsrpq$ is a (quasi-)Banach space which is continuously embedded in $\cD'$,
i.e. there exist $k$ and $c>0$ such that
$$
|\langle f, \phi\rangle| \le c\|f\|_{\Fsrpq}\cN_k(\phi)
\quad \mbox{for all} \;\; f\in\Fsrpq, \; \phi\in\cD.
$$

\medskip

We next introduce the sequence spaces $\fsrpq$ associated to $\Fsrpq$.
Here we assume that $\{\cX_j\}_{j=0}^\infty$ and $\cX:= \cup_{j=0}^\infty \cX_j$
are the sets of points from the definition of needles
with associated neighborhoods $\{I_\xi\}$, given in (\ref{def-I-xi}).


\begin{definition}
Suppose $s, \r \in \RR$, $0<p<\infty$, and $0<q\le\infty$. Then~$\fsrpq$
is defined as the space of all complex-valued sequences
$h:=\{h_{\xi}\}_{\xi\in \cX}$ such that
\begin{equation}\label{def-f-space}
\norm{h}_{\fsrpq} :=\nnorm{\Big(\sum_{j=0}^\infty
2^{sjq}\sum_{\xi \in \cX_j}
[|h_{\xi}|\WW(2^j;\xi)^{-\r/d}\tONE_{I_\xi}(\cdot)]^q\Big)^{1/q}}_{\Lp} <\infty
\end{equation}
with the usual modification for $q=\infty$. Here as before
$\tONE_{I_\xi}:=\mu(I_\xi)^{-1/2}\ONE_{I_\xi}$.
\end{definition}

The ``analysis" and ``synthesis" operators associated to the needlet systems $\Phi$, $\Psi$
are defined by
\begin{equation}\label{anal_synth_oprts}
S_\varphi: f\rightarrow \{\ip{f, \varphi_\xi}\}_{\xi \in \cX}
\quad\text{and}\quad
T_\psi: \{h_\xi\}_{\xi \in \cX}\rightarrow \sum_{\xi\in \X}h_\xi\psi_\xi.
\end{equation}
As in \cite{KPX1} one shows that the operator $T_\psi$ is well defined on $\fsrpq$,
namely, for any $h\in  \fsrpq$,
$T_\psi h:=\sum_{\xi\in \cX}h_\xi \psi_\xi$ converges in $\cD'$.
Moreover, the operator
$T_\psi: \fsrpq \to \cD'$ is continuous, i.e.
there exist constants $k>0$ and $c>0$ such that
\begin{equation}\label{synthesis}
|\langle T_\psi h, \phi\rangle| \le c \cN_k(\phi)\|h\|_{\fsrpq},
\quad \mbox{for all}\;\; h\in\fsrpq, \; \phi\in \cD.
\end{equation}

Our main result in this section asserts that the weighet $F$-spaces
can be characterized by the needlet coefficients of the distributions.


\begin{theorem}\label{thm:Fnorm-equivalence}
Let $s, \r\in \RR$, $0< p< \infty$ and $0<q\le \infty$.
The operators
$S_\varphi:\Fsrpq\rightarrow\fsrpq$ and $T_\psi:\fsrpq\rightarrow \Fsrpq$
are bounded and $T_\psi\circ S_\varphi=Id$ on $\Fsrpq$.
Consequently, $f\in \Fsrpq$ if and only if
$\{\ip{f, \varphi_\xi}\}_{\xi \in \cX}\in \fsrpq$.
Furthermore,
\begin{align}\label{Fnorm-equivalence-1}
\norm{f}_{\Fsrpq}
\sim  \norm{\{\ip{f,\varphi_\xi}\}}_{\fsrpq}
\sim \nnorm{\Big(\sum_{j=0}^\infty 2^{sjq}
\sum_{\xi\in \cX_j}[|\ip{f, \varphi_\xi}|\WW(2^j; \xi)^{-\r/d}|\psi_\xi(\cdot)|]^q\Big)^{1/q}}_{\Lp}.
\end{align}
In addition, the definition of $\Fsrpq$ is independent of the particular selection of
the type $(b)$ cutoff function $\hA$ satisfying $(\ref{star})$.
\end{theorem}

To us the spaces $\F sspq$ are more natural than the spaces
$\Fsrpq$ with $\r\ne s$ since they  embed ``correctly" with respect
to the smoothness index $s$.

%
\begin{proposition}\label{prop:F-embedding}
Let $0<p<p_1<\infty$, $0<q, q_1\le\infty$, and
$-\infty<s_1<s<\infty$. Then we have the continuous embedding
\begin{equation}\label{F-embed}
\F sspq \subset \F {s_1}{s_1}{p_1}{q_1} \quad\mbox{if}\quad
s/d-1/p=s_1/d-1/p_1.
\end{equation}
\end{proposition}

The proof of this proposition is quite similar to the proof of the respective
embedding result on $B^d$ in \cite[Proposition 4.11]{KPX2} and will be omitted.

\medskip

We have the following identification of spaces $\F 00p2$.


\begin{proposition}\label{prop:ident}
We have
$$
\F 00p2\sim \Lpp,
\quad 1 < p < \infty,
$$
with equivalent norms.
Consequently, for any $f\in\Lpp$, $1<p<\infty$,
$$
\|f\|_p
\sim \nnorm{\Big(\sum_{j=0}^\infty
\sum_{\xi\in \cX_j}(|\ip{f, \varphi_\xi}||\psi_\xi(\cdot)|)^2\Big)^{1/2}}_{\Lp}.
$$
\end{proposition}

The proof of this proposition uses the multipliers from Theorem~\ref{thm:multipliers}
and can be carried out exactly as in the case of spherical
harmonic expansions in \cite[Proposition~4.3]{NPW2}.
We omit~it.

\section{Weighted Besov spaces on $[-1, 1]^d$}\label{Besov}
\setcounter{equation}{0}

To define weighted Besov spaces on $[-1, 1]^d$ we use again
the sequence of kernels $\{\Phi_j\}$ introduced in (\ref{def-Phi-j}) with
$\hA$ a cutoff function of type (b) obeying (\ref{star}).
We shall keep the development of these spaces short since
the proofs of the results are the same as in the univariate case,
given in \cite{KPX1}.


\begin{definition}\label{def-Besov}
Let $s, \r\in \RR$ and $0<p,q \le \infty$. The weighted Besov space
$\Bsrpq := \Bsrpq(w_{\a,\b})$ is defined
as the set of all $f \in \cD'$ such that
\begin{equation}\label{def-Besov-sp}
\|f\|_{\Bsrpq} :=
\Big(\sum_{j=0}^\infty \Big(2^{s j}
\|\WW(2^j; \cdot)^{-\r/d}\Phi_j*f(\cdot)\|_{\Lp}\Big)^q\Big)^{1/q}
< \infty,
\end{equation}
where the $\ell^q$-norm is replaced by the sup-norm if $q=\infty$.
\end{definition}

Note that as in the case of weighted Triebel-Lizorkin spaces the above definition
is independent of the particular choice of $\hA$ and  $\Bsrpq$ is a (quasi-)Banach
space which is continuously embedded in $\cD'$.

\smallskip

We next introduce the sequence spaces $\bsrpq$ associated to $\Bsrpq$.
To this end we use some of the notation established in the previous section.


\begin{definition}
Let $s, \r\in \RR$ and $0<p,q \le \infty$. Then $\bsrpq$
is defined to be the space of all complex-valued sequences
$h:=\{h_{\xi}\}_{\xi\in \cX}$ such that
\begin{equation}\label{def-berpq}
\norm{h}_{\bsrpq} :=\Bigl(\sum_{j=0}^\infty
2^{j(s-d/p+d/2)q}
\Bigl[\sum_{\xi\in \cX_j}\Big(\WW(2^j;\xi)^{-\r/d+1/p-1/2}|h_\xi|\Big)^p
\Bigr]^{q/p}\Bigr)^{1/q}
\end{equation}
is finite,
with the usual modification for $p=\infty$ or $q=\infty$.
\end{definition}

The analysis and synthesis operators
$S_\ph$ and $T_\psi$ defined in (\ref{anal_synth_oprts}) play an important role here.
As for weighted Triebel-Lizorkin spaces the operator $T_\psi$ is well defined on $\bsrpq$,
i.e. for any $h\in  \bsrpq$,
$T_\psi h:=\sum_{\xi\in \cX}h_\xi \psi_\xi$ converges in $\cD'$.
Also, the operator $T_\psi: \bsrpq \to \cD'$ is continuous.

The following characterization of weighted Besov spaces is the main result of this section.


\begin{theorem}\label{thm:Bnorm-eq}
Let $s, \r\in \RR$ and  $0< p,q\le \infty$.
Then the operators
$S_\varphi:\Bsrpq\rightarrow\bsrpq$ and
$T_\psi:\bsrpq\rightarrow \Bsrpq$
are bounded and $T_\psi\circ S_\varphi=Id$ on $\Bsrpq$.
Consequently, for $f\in \cD'$ we have that $f\in \Bsrpq$ if and
only if $\{\ip{f, \varphi_\xi}\}_{\xi \in \cX}\in \bsrpq$.
Moreover,
$$
\| f\|_{\Bsrpq}  \sim  \norm{\{\ip{f,\varphi_\xi}\}}_{\bsrpq}
\sim \Big(\sum_{j=0}^\infty2^{sjq} \Bigl[\sum_{\xi\in \cX_j}
\Big(\WW(2^j; \xi)^{-\r/d}\norm{\ip{f,\varphi_\xi}\psi_\xi}_{\Lp}\Big)^p\Bigr]^{q/p}\Bigr)^{1/q}.
$$
In addition, the definition of $\Bsrpq$ is independent of the particular selection of
the type $(b)$ cutoff function $\hA$ satisfying $(\ref{star})$.
\end{theorem}

The parameter $\r$ in the definition of $\Bsrpq$ allows to consider various scales of weighted Besov spaces.
The spaces $\B s0pq$ can be regarded as ``classical" Besov spaces.
However, to us more natural are the spaces $\B sspq$ ($\r=s$) which in contrast to $\B s0pq$,
first, embed ``correctly" with respect to the smoothness index $s$,
and secondly, the right smoothness spaces in nonlinear $n$-term weighted approximation from needles
are defined in terms of spaces $\B sspq$ (see \S\ref{Nonlin-app} below).

%
\begin{proposition}\label{B-embedding}
Let $0<p\le p_1\le\infty$, $0<q\le q_1\le \infty$, and
$-\infty<s_1\le s<\infty$. Then we have the continuous embedding
\begin{equation}\label{B-embed}
\B sspq \subset \B {s_1}{s_1}{p_1}{q_1} \quad\mbox{if}\quad
s/d-1/p=s_1/d-1/p_1.
\end{equation}
\end{proposition}

This proposition is an immediate consequence of estimate (\ref{norm-relation2}).

\section{Application of weighted Besov spaces to nonlinear approximation}
\label{Nonlin-app}
\setcounter{equation}{0}

We now consider nonlinear $n$-term approximation for a needlet system
$\{\psi_\eta\}_{\eta\in \cX}$ with $\ph_\eta=\psi_\eta$,
defined as in (\ref{def.Phi-j})-(\ref{def-needlets2})
with $\hB = \hA$, $\hA\ge 0$,
i.e. $\hA\ge 0$ is a first or second kind admissible cutoff function  of type (c)
(see Definitions~\ref{cutoff-first}-\ref{cutoff-second}).
Then $\{\psi_\eta\}$ are real-valued.

Let $\Sigma_n$ be the nonlinear set of all functions $g$ of the form
$$
g = \sum_{\xi \in \Lambda} a_\xi \psi_\xi,
$$
where $\Lambda \subset \cX$, $\#\Lambda \le n$,
and $\Lambda$ is allowed to vary with $g$.
Denote by $\sigma_n(f)_p$ the error of best $\Lpp$-approximation to
$f \in \Lpp$ from $\Sigma_n$, i.e.
$$
 \sigma_n(f)_p := \inf_{g \in \Sigma_n} \|f - g\|_p.
$$
We consider approximation in $\Lpp$, $0<p<\infty$.

Assume $0 < p < \infty$, $s > 0$, and $1/\tau := s/d + {1}/{p}$
and denote briefly $B^s_\tau:= B^{ss}_{\tau\tau}$.
By Theorem~\ref{thm:Bnorm-eq} and (\ref{norm-Needlets})
it follows that
\begin{equation}\label{Btau-norm}
\|f\|_{B^s_\tau}\sim
\Big(\sum_{\xi\in\cX} \|\langle f, \psi_\xi \rangle \psi_\xi\|_p^\tau\Big)^{1/\tau}.
\end{equation}
Exactly as in \cite[Proposition 6.1]{NPW2} this leads to
the embedding of $B^s_\tau$ into $\Lpp$,
which plays an important role in the proof of the main result of this section:


\begin{theorem}\label{thm:jackson} {\bf [Jackson estimate]}
If $f \in B^s_\tau$, then
\begin{equation}\label{jackson}
\sigma_n(f)_p \le c n^{-s/d}\|f\|_{B^s_\tau}, \quad n\ge 1,
\end{equation}
where $c>0$ depends only on $s$, $p$,
and $\hA$.
\end{theorem}

The proofs of this theorem can be carried out exactly as the proofs
of the Jackson estimate in \cite[Theorem 6.2]{NPW2}. We omit it.

It is an important open problem to prove the companion to (\ref{jackson})
Bernstein estimate: If $g \in \Sigma_n$ and $1<p<\infty$, then
\begin{equation}\label{bernstein}
\|g\|_{B^s_\tau} \le c n^{s/d} \|g\|_p.
\end{equation}
If true this estimate would enable one to characterize the rates (approximation spaces)
of nonlinear $n$-term approximation in $\Lpp$ ($1<p<\infty$) from needlet systems.

\section{Weighted Triebel-Lizorkin and Besov spaces on $B^{d_1}\times [-1,1]^{d_2}$}
\label{cross-products}
\setcounter{equation}{0}

Our aim is to briefly describe how the theory of weighted spaces of distributions on the product set
$B^{d_1}\times [-1,1]^{d_2}$
can be developed via tensor product orthogonal polynomials.

\subsection{Localized kernels for orthogonal polynomials on the ball}
Localized polynomial kernels on the unit ball $B^d$ in $\RR^d$ have been developed in \cite{PX2}
and utilized in \cite{KPX2} to the development of Triebel-Lizorkin and Besov spaces
on $B^d$ with weight
$$
w_\mu(x):=(1-\|x\|_2^2)^{\mu-1/2}, \quad \mu \ge 0.
$$
Here, we compile all needed results from \cite{KPX2, PX2} and give some new facts.

Denote by $V_n$
the set of all polynomials of degree $n$ in $d$ variables which are orthogonal to
the lower degree polynomials in $L^2(B^{d}, w_\mu)$ and
let $P_n(w_\mu, x, y)$ be the kernel
of the orthogonal projector $\Proj_n: L^2(B^{d}, w_\mu) \to V_n$, i.e.
\begin{equation}\label{def-Pn}
(\Proj_n f)(x) = \int_{B^d} f(y)P_n(w_\mu; x, y) w_\mu(y) dy.
\end{equation}
An explicit representation of the reproducing kernel $P_n(w_\mu, x, y)$ is given
in \cite{X99}:
\begin{equation} \label{ball-reprod-kernel}
P_n(w_\mu;x,y) = c_\mu \frac{n+\lambda}{\lambda} \int_{-1}^1
     C_n^\lambda (t(x,y;u)) (1-u^2)^{\mu-1} du,
\quad\mbox{$\; \mu > 0$,}
\end{equation}
where
$$
 t(x,y;u):= \langle x, y \rangle + u\sqrt{1-\|x\|_2^2}\sqrt{1-\|y\|_2^2} \quad \hbox{and}
  \quad  \lambda: = \mu + \frac{d-1}{2}.
$$
For an admissible univariate cutoff function $\ha$ (see Definition~\ref{cutoff-d1}), denote
\begin{equation} \label{ball-kernel}
L_n^\mu(x, y):= \sum_{j=0}^\infty \ha\Big(\frac{j}{n}\Big) P_j(w_\mu;x,y)
\end{equation}
Analogues of Theorems~\ref{thm:product-Jacobi-0} and \ref{thm:product-Jacobi-1} on $B^d$
are established in \cite{PX2}.
Denote
\begin{equation}\label{ball-A-W}
W_\mu(n;x) := \Big(\sqrt{1-\|x\|_2^2} + n^{-1}\Big)^{2 \mu}
\end{equation}
and
\begin{equation}\label{def-dist-ball}
\rho(x,y):
=\arccos \Big(\langle x, y \rangle + \sqrt{1-\|x\|_2^2}\sqrt{1-\|y\|_2^2}\Big),
\end{equation}
which is a distance on $B^d$.


\begin{theorem}\label{thm:local-ball}
Given an admissible univariate function $\ha$, 
for any $\sigma>0$ there exists a constant $c>0$ such that
\begin{equation} \label{est-ball-kernel}
|L_n^\mu(x,y)|\le\frac{c n^d }{\sqrt{W_\mu(n; x)}\sqrt{W_\mu(n; y)}}
\big(1 + n\rho(x, y)\big)^{-\sigma},
\quad x, y\in B^d.
\end{equation}
Furthermore, for any $x, y, \xi\in B^d$ such that $\rho(x, \xi) \le c_*n^{-1}$
\begin{equation} \label{ball-kernel-Lip}
|L_n^\mu(x,y)-L_n^\mu(\xi,y)|\le\frac{c n^{d+1}\rho(x, \xi) }{\sqrt{W_\mu(n; x)}\sqrt{W_\mu(n; y)}}
\big(1 + n\rho(x, y)\big)^{-\sigma}.
\end{equation}
\end{theorem}

This theorem was established in \cite{PX2} (Theorem 4.2 and Proposition~4.7)
in the case of admissible cutoff functions $\ha$
which are constant around $t=0$. Its proof hinges on the localization of the kernels $\cQ_n^{\a,\b}$
from \eqref{eq:L-n}. Due to Theorem~\ref{thm:est-Q} now Theorem~\ref{thm:local-ball} holds
for admissible cutoff functions $\ha$ in the sense of Definition~\ref{cutoff-d1}
with the proof from \cite{PX2}.

We shall need two additional estimates with the first being the analogue of
Lemma~\ref{lem:simple-est} on $B^d$.


\begin{lemma} \label{ball-lem1}
For $x,y \in B^d$,
\begin{equation} \label{ball-Pn}
|P_n(w_\mu;x,y)| \le  \frac{cn^{d-1}}{\sqrt{W_\mu(n;x)} \sqrt{W_\mu(n;y)}}.
\end{equation}
\end{lemma}

\noindent
{\bf Proof.}
The proof of this lemma relies on the following estimate that follows
from Theorem 3.1 in \cite{DaiXu}:
If $a$ and $b$ are constants such that $|a|+|b| \le 1$, then
$$
   \left| \int_{-1}^1 C_n^\lambda (a u + b) (1-u^2)^{\mu-1} du  \right|
     \le c n^{2 \lambda - 2 \mu -1} \frac{\big(|a| + n^{-1} \sqrt{1-|a| -|b|} + n^{-2}\big)^{-\mu}}
         {\big(1 + n \sqrt{1-|a|-|b|}\big)^{\lambda -\mu}}.
$$
Denote briefly $A(x): = \sqrt{1-\|x\|_2^2}$.
We apply the above inequality with $a = A(x)A(y)$ and $b = \langle x,y\rangle$.
Setting $\|x\|_2 =: \cos \theta$ and $\|y\|_2 =: \cos \phi$, $0 \le \theta, \phi \le \pi$,
we have
\begin{align*}
  1- |a|-|b| & \ge 1 - \|x\|_2 \cdot \|y\|_2 - \sqrt{1-\|x\|_2^2}\sqrt{1-\|y\|_2^2} =
      1 - \cos (\theta - \phi) \\
        & = 2 \sin^2 \frac{\theta -\phi}2 \ge c (\theta - \phi)^2 \ge c (\sin \theta - \sin \phi)^2
         = c \left(A(x) - A(y)\right)^2,
\end{align*}
and hence
\begin{equation} \label{ball-Pn2}
 | P_n(w_\mu;x,y)| \le c n^{2 \lambda - 2 \mu} \left(A(x)A(y) + n^{-1}
       |A(x) - A(y)|  + n^{-2}\right)^{-\mu}.
\end{equation}
Here we used that $\big(1+ n \sqrt{1-|a| -|b|}\big)^{\lambda - \mu} \ge 1$.
Now, from $A(x), A(y)\ge 0$ it easily follows that
\begin{equation} \label{ball-1}
A(x)A(y) + n^{-1} |A(x) - A(y)|  + n^{-2}  \sim ( A(x) + n^{-1})( A(y) + n^{-1}).
\end{equation}
This coupled with \eqref{ball-Pn2} yields \eqref{ball-Pn}.
$\qed$

\smallskip

The next lemma gives an analogue of estimate \eqref{Pn-Pn} on the ball.


\begin{lemma} \label{ball-lem3}
For any $x,y,\xi \in B^d$ such that $\rho (x,\xi) \le c^* n^{-1}$,
\begin{equation}\label{ball-Lip}
   |P_n(w_\mu;x,y) - P_n(w_\mu;\xi,y)| \le  \frac{c n^{d} \rho(x,\xi)}
        {\sqrt{W_\mu(n;x)} \sqrt{W_\mu(n;y)}},
\end{equation}
where the constant $c>0$ depends only on $\mu$, $d$, and $c^*$.
\end{lemma}

The proof of this lemma is somewhat lengthy and will be given in the appendix.

\subsection{Localized cross product basis kernels}

We consider orthogonal polynomials on
$B^{d_1}\times [-1,1]^{d_2}$ with weight
$$
w_{\mu,\a,\b}(x):=w_\mu(x')w_{\a,\b}(x''),
\quad x=(x', x''), \quad x'\in B^{d_1}, x''\in [-1,1]^{d_2},
$$
where $w_\mu(x'):=(1-\|x'\|_2^2)^{\mu-1/2}$, $\mu \ge 0$, and
$w_{\a, \b}(x''):=\prod_{j=1}^{d_2} w_{\a_j, \b_j}(x_j'')$
with $\a_j, \b_j\ge -1/2$ as in (\ref{Jacobi-weight}).

Denote by $V_n$ the set of all algebraic polynomials of degree $n$ in $d_1$ variables
which are orthogonal to the lower degree polynomials in
$L^2(B^{d_1}, w_\mu)$ and let $P_n(w_\mu, x', y')$ be the kernel
of the orthogonal projector
$\Proj_n: L^2(B^{d_1}, w_\mu) \to V_n$, see \eqref{def-Pn}-\eqref{ball-reprod-kernel}.

We are interested in kernels of the form
\begin{equation}\label{def-Lambda-n}
\Lambda_n(x, y)
:=\sum_{(j, \nu)\in \NN_0\times \NN_0^{d_2}}\hA\Big(\frac{j}{n},\frac{\nu}{n}\Big)
P_j(w_\mu; x', y')\tP_\nu^{(\a, \b)}(x'')\tP_\nu^{(\a, \b)}(y'').
\end{equation}
Here $\hA\in C^\infty[0, \infty)^{1+d_2}$ is an admissible cutoff function
in the sense of Definition~\ref{cutoff-first}
and $\tP_\nu^{(\a, \b)}$ are the tensor product Jacobi polynomials
defined as in (\ref{product-Jacobi}).
To estimate the localization of $\Lambda_n(x, y)$ we need the weight
\begin{equation}\label{modified-weight}
W_{\mu,\a,\b}(n;x):=W_\mu(n;x')\WW(n; x''),
\end{equation}
where $W_\mu(n;x')$ is defined as in \eqref{ball-A-W} and $\WW(n; x'')$
as in (\ref{Jacobi-weight-d}).
We~also need the distance $\rho_*(x, y)$ on $B^{d_1}\times [-1,1]^{d_2}$
defined by
\begin{align}\label{cross-distance}
\rho_*(x, y)
&:= \max\Big\{\rho(x', y'),
\max_{1\le j\le d_2}|\arccos x_j''-\arccos y_j''|\Big\},\notag
\end{align}
where $\rho(x', y')$ is the distance on $B^{d_1}$ defined as in \eqref{def-dist-ball}.

We now give the localization of the kernels $\Lambda_n(x, y)$ from (\ref{def-Lambda-n}):


\begin{theorem}\label{thm:product-ball-Jacobi}
If $\hA\in C^\infty[0, \infty)^{1+d_2}$ is an admissible cutoff function in the sense of Definition~\ref{cutoff-first},
then for any $\sigma>0$ there exists a constant $c>0$ such that
\begin{equation}\label{ball-Jacobi}
|\Lambda_n(x,y)|\le\frac{c n^{d_1+d_2} }{\sqrt{W_{\mu,\a,\b}(n; x)}\sqrt{W_{\mu,\a,\b}(n; y)}}
\big(1 + n\rho_*(x, y)\big)^{-\sigma}
\end{equation}
 for $x, y\in B^{d_1}\times[-1,1]^{d_2}$.
\end{theorem}

This theorem is an immediate consequence
of Theorems~\ref{thm:product-Jacobi-0},\ref{thm:local-ball}, and Lemma~\ref{ball-lem1}.
 (see the proof of Theorem~\ref{thm:product-Jacobi-0}).

\smallskip

The analogue of Theorem~\ref{thm:product-Jacobi-1} reads as follows:


\begin{theorem}\label{thm:product-ball-Jacobi-1}
Let $\hA$ be an admissible cutoff function which belongs to the class
$\cS(d_2+1,\cL;\gamma,\bargamma M)$
for some $\cL$ and $M$ as in \eqref{def-M} and $\gamma,\bargamma>0$
$($see Definition~\ref{Def-S}$)$.
Then the kernels from \eqref{def-Lambda-n} satisfy
\begin{equation} \label{ball-Jacobi-kernel}
|\Lambda_n(x, y)|
\le \frac{c n^{d_1+d_2} }{\sqrt{W_{\mu,\a,\b}(n; x)}\sqrt{W_{\mu,\a,\b}(n; y)}}
\exp\Big\{-\frac{\tilde{c}n\rho_*(x, y)}{\cL(n\rho_*(x, y))}\Big\}
\end{equation}
for $x, y\in B^{d_1}\times[-1,1]^{d_2}$.
Here $\tilde{c}=c'/\bargamma M$, where $c'>0$ is an absolute constant.
\end{theorem}
Here the argument is the same as for the proof of Theorem~\ref{thm:product-Jacobi-1}
(see the proof of Theorem~\ref{thm:product-Jacobi-0}) using Theorem~\ref{thm:product-Jacobi-1}
and \cite[Theorem~6.1]{IPX}.

An analogue of Theorem~\ref{thm:Lip} is also valid:


\begin{theorem}\label{thm:Lip-cross}
If $\hA\in C^\infty[0, \infty)^{1+d_2}$ is an admissible cutoff function, then
for any $\sigma>0$ and
for all $x, y, \xi\in B^{d_1}\times[-1,1]^{d_2}$
such that $\rho(x, \xi) \le c_*n^{-1}$, $n\ge 1$, $c_*>0$,
the kernel $\Lambda_n$ from $(\ref{def-Lambda-n})$ satisfies
\begin{equation}\label{Lip-cross}
|\Lambda_n(x,y)-\Lambda_n(\xi,y)|
\le\frac{c n^{d_1+d_2+1}\rho(x, \xi)}{\sqrt{\WWW(n; x)}\sqrt{\WWW(n; y)}}
\big(1+ n\rho(x, y)\big)^{-\sigma},
\end{equation}
where $c>0$ depends only on $\sigma, d, \a, \b$, $c_*$, and $\hA$.
\end{theorem}

The proof of this theorem is quite similar to the proof of Theorem~\ref{thm:Lip} and
relies on Theorems~\ref{thm:Lip},\ref{thm:local-ball}, and Lemma~\ref{ball-lem3}.

\subsection{Construction of needlets on \boldmath $B^{d_1}\times [-1,1]^{d_2}$}

An important component of our theory is the construction of frames on $B^{d_1}\times [-1,1]^{d_2}$.
To this end one uses a Calde\'{o}n type formula based on localized kernels as the kernels
in \eqref{def-Lambda-n} and a cubature formula.
A cubature formula on $B^{d_1}\times [-1,1]^{d_2}$ exact for sufficiently large degree polynomials
can be constructed as product of the cubature formula on $B^{d_1}$ from \cite[\S5]{PX2}
and the cubature on $[-1,1]^{d_2}$ from \S\ref{sec:cubature}.
Once the components are in place, the construction is carried out exactly as in \S\ref{sec:Needlets}.
We skip the details.

\subsection{Spaces of distributions on \boldmath $B^{d_1}\times [-1,1]^{d_2}$}

It is natural to use as test functions the set $\cD:=C^\infty(B^{d_1}\times [-1,1]^{d_2})$,
where the topology is defined by the semi-norms
$
|\phi|_\mu:=\|D^\mu \phi\|_\infty
$
for all multi-indices $\mu$.
Just as in the case of tensor product Jacobi polynomials (\S\ref{distributions})
the test functions $\phi\in\cD$ can be characterized by their cross polynomial expansions on
$B^{d_1}\times [-1,1]^{d_2}$.
The space $\cD'$ of distributions on $B^{d_1}\times [-1,1]^{d_2}$ is defined as the set
of all continuous linear functionals on $\cD$.

For an admissible cutoff function $\hA:[0, \infty)^{1+d_2}\mapsto \bC$ of type (b)
obeying condition (\ref{star}) we define 
$\Phi_0(x, y) := P_0(w_\mu; x', y')\tP_0^{(\a, \b)}(x'')\tP_0^{(\a, \b)}(y'')$ and
$$
\Phi_j(x, y)
:=\sum_{(m, \nu)\in \NN_0\times \NN_0^{d_2}}\hA\Big(\frac{m}{2^{j-1}},\frac{\nu}{2^{j-1}}\Big)
P_m(w_\mu; x', y')\tP_\nu^{(\a, \b)}(x'')\tP_\nu^{(\a, \b)}(y''),
\; j\ge 1.
$$

Then the weighted Triebel-Lizorkin space $\Fsrpq:=\Fsrpq(w_{\mu,\a,\b})$
with  $s, \r \in \RR$, $0<p<\infty$, and $0<q\le\infty$,
is defined as the set of all $f\in \cD'$ such that
\begin{equation}\label{Tri-Liz-norm-cross}
\|f\|_{\Fsrpq}:=\Big\|\Big(\sum_{j=0}^{\infty}
\Big[2^{sj}\WWW(2^j;\cdot)^{-\r/(d_1+d_2)}|\Phi_j*f(\cdot)|\Big]^q\Big)^{1/q}\Big\|_{\Lp} <\infty
\end{equation}
with the usual modification when $q=\infty$.
Here $\Phi_j*f$ is defined as in \eqref{convolution1}.

The weighted Besov space $\Bsrpq := \Bsrpq(w_{\mu,\a,\b})$
with $s, \r\in \RR$ and $0<p,q \le \infty$,
is defined as the set of all $f \in \cD'$ such that
\begin{equation}\label{Besov-norm-cross}
\|f\|_{\Bsrpq} :=
\Big(\sum_{j=0}^\infty \Big[2^{s j}
\|\WWW(2^j; \cdot)^{-\r/(d_1+d_2)}\Phi_j*f(\cdot)\|_{\Lp}\Big]^q\Big)^{1/q}
< \infty,
\end{equation}
where the $\ell^q$-norm is replaced by the sup-norm if $q=\infty$.

Without going into further details, we note that the theory of
Triebel-Lizorkin and Besov space on $B^{d_1}\times [-1,1]^{d_2}$ with weight $w_{\mu,a,b}(x)$
can be further developed in analogy to the spaces on $[-1,1]^d$ from \S\S\ref{Tri-Liz}-\ref{Besov}.
Also, needlets on $B^{d_1}\times [-1,1]^{d_2}$ can be deployed for the decomposition of
the $F$- and $B$-spaces on $B^{d_1}\times [-1,1]^{d_2}$ as in \S\S\ref{Tri-Liz}-\ref{Besov}.
The point is that all ingredients needed for this theory are either in place
or can easily be developed.

\section{Discussion}\label{discussion}
\setcounter{equation}{0}

Although this paper is mainly concerned with weighted Triebel-Lizorkin and Besov space on $[-1,1]^d$
it is one of our goals to show how the theory of $F$- and $B$-spaces can be developed on products of
$[-1,1]^{d_1}$, $B^{d_2}$, $\SS^{d_3}$, $T^{d_4}$, $\RR^{d_5}$, or $\RR^{d_6}_+$ with weights.
For $B^{d_1}\times [-1,1]^{d_2}$ a sketch of the main ingredients of the theory was given
in the previous section.
We belive that the most natural way to define and develop this sort of spaces is
via orthogonal decompositions, where kernels like the ones from \eqref{def.L},
\eqref{Jacobi-kernel-d} or \eqref{def-Lambda-n} play a prominent role.

We would like to turn again our attention to the fundamental question of what kind of
cutoff functions $\hA$ can be used in the case of cross product bases.
As was already mentioned in the introduction, as for univariate Jacobi polynomials (see \eqref{def.L})
univariate cutoff functions $\ha$ induce rapidly decaying kernels on
the sphere \cite{NPW2}, ball \cite{KPX2}, simplex \cite{IPX},
and in the context of tensor product Hermite \cite{PX3} and Laguerre functions \cite{KPPX}.
Note that cutoff functions $\ha$ which are constants around $t=0$ are sufficient
for the development of the theory in these cases.
However, as was already seen truly multivariate cutoff functions $\hA$ need to be used in
the case of product Jacobi polynomials or cross product bases.
Moreover, the localization of the respective kernels depends on the behavior of $\hA$
at the boundary of $[0, \infty)^d$, i.e. at the coordinate planes.
This is intimately related to the impact of the behavior of the univariate cutoff functions
$\ha$ at $t=0$ on the localization of the kernels on the interval, ball, sphere, etc.
This behavior appears as a boundary condition on $\hA$ and becomes an important issue.

The key observation is that (as in Theorem~\ref{thm:Jacobi-localization-0})
the localization results given in the theorems described below hold under the condition that
the compactly supported $C^\infty$ univariate cutoff function $\ha$ satisfies
$$
\ha^{(m)}(0)=0 \quad \mbox{for $\;m=1, 2, \dots$}
$$
These are:
(1) Theorem~4.2 in \cite{PX2} on the ball,
(2) Theorem~2.2 in \cite{NPW2} on the sphere,
(3) Theorem~7.1 in \cite{IPX} on the simplex,
(4) Corollary~1 in \cite{PX3} for tensor product Hermite functions,
(5) Theorems~3.2,~3.7,~3.8 in \cite{KPPX} for tensor product Laguerre functions.
The proofs of these results utilize the scheme of the proof of Theorem~\ref{thm:Jacobi-localization-0}
with very little variations and will be omitted.
Consequently, the cross product basis kernels induced by an admissible cutoff function $\hA$
(see Definition~\ref{cutoff-first}) obtained from any combination of the above mentioned bases on
$[-1,1]^{d_1}$, $B^{d_2}$, $\SS^{d_3}$, $T^{d_4}$, $\RR^{d_5}$, or $\RR^{d_6}_+$
will decay rapidly as in Theorems~\ref{thm:product-Jacobi-0}, \ref{thm:product-ball-Jacobi}.
Further modifications and extensions as in
Theorems~\ref{thm:product-Jacobi-1}, \ref{thm:Lip}, \ref{thm:product-ball-Jacobi-1},
\ref{thm:Lip-cross} are also almost automatic.

The construction of needlets on products of two or more of the sets
$[-1,1]^{d_1}$, $B^{d_2}$, $\SS^{d_3}$, $T^{d_4}$, $\RR^{d_5}$, or $\RR^{d_6}_+$
follows easily the pattern of the construction on $[-1, 1]^d$ from \S\ref{Needlets},
based on tensor product basis kernels and product cubature formulas.

The ensuing program for developing weighted Triebel-Lizorkin and Besov spaces
on products of sets as above can be carried out as for the spaces on $[-1, 1]^d$
developed in this article.

\section{Appendix: Proof of Lemma~\ref{ball-lem3}.}\label{appendix}
For $\mu =0$ the expression of $P_n(w_\mu;x,y)$ in \eqref{ball-reprod-kernel}
simplifies considerably as $\mu \to 0$;
the integral becomes a sum of two terms, as shown in \cite{X99}.
This case is easier than the case $\mu >0$. We omit its proof.

Assume $\mu  > 0$.
The proof hinges on the following lemma which is an immediate consequence of Lemma 3.5 in \cite{DaiXu}.

\begin{lemma} \label{ball-lem2}
Suppose $\mu > 0$, $0< |a| < 1$, $\eta\in C^\infty [-1,1]$ with
$\supp \eta \subset [-\frac 12,1]$.
If $|b| \le 1 - |a|$, then
$$
  \left|\int_{-1}^1 C_n^\lambda (a t + b) \eta(t) (1-t)^{\mu -1} dt \right|
  \le  \frac{cn^{2 \lambda - 2 \mu -1} }{|a|^{\mu}
              \left(1+n \sqrt{1-|a+b|}\right)^{\lambda - \mu}}.
$$
\end{lemma}

The proof of Lemma~\ref{ball-lem3} will be divided into two parts.

\medskip
{\it Case 1:} $A(x)A(y) \le  16c^*(n^{-1} \rho(x,y) + n^{-2})$,
where $c^*$ is the constant from the hypothesis of Lemma~\ref{ball-lem3}.
We shall need the following estimate for Gegenbauer polynomials,
which follows from \eqref{est.Pn2}:
\begin{equation} \label{ball-Gegen}
|C_n^\lambda(t)| \le c n^{2\lambda -1} \left(1+ n \sqrt{1-t^2} \right)^{- \lambda},
\quad t \in [-1,1].
\end{equation}
Denote by $I_u$ the interval with end points $t(x,y;u)$ and $t(\xi,y;u)$.
Then using the identity
$\frac{d}{dt} C_n^\lambda(t) = 2\lambda C_{n-1}^{\lambda+1}(t)$ \cite[(4.7.27)]{Sz}, we obtain
\begin{align} \label{ball-3}
 E: &= | P_n(w_\mu;x,y) - P_n(w_\mu;\xi,y)| \notag\\
 & \le c n \int_{-1}^1 \left |C_n^\lambda(t(x,y;u)) -C_n^\lambda(t(\xi,y;u))\right|
       (1-u^2)^{\mu-1} du \\
 & \le  c n \int_{-1}^1  \|  C_{n-1}^{\lambda+1}\|_{L^\infty(I_u)} |t(x,y;u) - t(\xi,y;u)|
       (1-u^2)^{\mu-1} du. \notag
\end{align}
By \eqref{ball-Gegen} it follows that
\begin{align}\label{ball-44}
    \|  C_{n-1}^{\lambda+1}\|_{L^\infty(I_u)} \le  c \,n^{2\lambda +1}
        & \left[ \left(1+ n \sqrt{1-t(x,y;u)^2} \right)^{- \lambda-1}  \right. \\
           & \qquad  + \left.
               \left(1+ n \sqrt{1-t(\xi,y;u)^2} \right)^{- \lambda-1} \right] \notag.
\end{align}
%
%
%
%
If $t(x,y;u)\ge 0$, then
\begin{align*}
1-t(x,y;u)^2
&\ge 1-t(x,y;u) \ge 1 - \langle x, y \rangle - A(x)A(y)\\
& = 1- \cos\rho(x, y) = 2\sin^2(\rho(x, y)/2)\ge (2/\pi^2)\rho (x, y)^2,
\end{align*}
and similarly if $t(x,y;u)<0$, then
\begin{align*}
1-t(x,y;u)^2
&\ge 1+t(x,y;u) \ge 1 + \langle x, y \rangle - A(x)A(y)\\
&= 1 - \langle x, -y \rangle - A(x)A(-y) \ge (2/\pi^2)\rho (x, -y)^2.
\end{align*}
The above estimates along with \eqref{ball-44} and $\rho(x,\xi) \le c^* n^{-1}$ yield
\begin{align*}
&\|C_{n-1}^{\lambda+1}\|_{L^\infty(I_u)} \le  c \,n^{2\lambda +1}
\left(1+ n \rho(x,y) \right)^{- \lambda-1}
\quad \mbox{if}\quad t(x,y;u)\ge 0, \quad \mbox{and}\\ 
&\|C_{n-1}^{\lambda+1}\|_{L^\infty(I_u)} \le  c \,n^{2\lambda +1}
\left(1+ n \rho(x,-y) \right)^{- \lambda-1}
\quad \mbox{if}\quad t(x,y;u)< 0. 
\end{align*}
We use these inequality in \eqref{ball-3} to obtain
\begin{align*}
E \le  cn^{2\lambda+2}
&\int_{-1}^1 \frac{|t(x,y;u) - t(\xi,y;u)|}
                   {\left(1+ n \rho(x,y) \right)^{ \lambda+1} } (1-u^2)^{\mu-1} du\\
+  c n^{2\lambda+2}
&\int_{-1}^1 \frac{|t(x,y;u) - t(\xi,y;u)|}
    { \left(1+ n \rho(x, -y) \right)^{ \lambda+1} } (1-u^2)^{\mu-1} du
=: E_1+E_2.
\end{align*}

To estimate $E_1$ and $E_2$ we shall need the inequality (see \cite[Lemma 4.1]{PX2})
\begin{equation}\label{Ax-Ay}
|A(x)-A(y)| \le \sqrt{2}\rho(x, y),
\quad x, y\in B^d,
\end{equation}
which implies
\begin{align}\label{Ax<Ay-rho}
A(y)+n^{-1}
&\le A(x)+n^{-1}+ \sqrt{2}\rho(x, y)\\
&\le \sqrt{2}(A(x)+n^{-1})(1+n\rho(x, y)).\notag
\end{align}
On the other hand, by \eqref{ball-1}, \eqref{Ax-Ay}, and our assumption
it follows that
\begin{align}\label{AxAy}
(A(x)+n^{-1})(A(y)+n^{-1})
&\le c\big(A(x)A(y)+ n^{-1}|A(x)-A(y)|+n^{-2}\big)\\
&\le cn^{-2}(1+n\rho(x, y)).\notag
\end{align}
This along with \eqref{Ax<Ay-rho} gives
\begin{equation}\label{Ax<rho}
A(y)^2 \le \sqrt{2}(A(x)+n^{-1})(A(y)+n^{-1})(1+n\rho(x, y))
\le cn^{-2}(1+n\rho(x, y))^2.
\end{equation}
As in \cite[p. 136]{PX2} we have using \eqref{Ax-Ay}
\begin{align}\label{t-t}
|t(x,y;u) - t(\xi,y;u)|
&\le |\cos\rho(x, y)- \cos\rho(\xi, y)| + \sqrt{2}|1-u|A(y)\rho(x, \xi)\notag\\
&\le \rho(x, \xi)\big(\rho(x, y)+\rho(\xi, y)\big) + \sqrt{2}A(y)\rho(x, \xi).
\end{align}
Combining this with \eqref{Ax<rho} and $\rho(x,\xi) \le c n^{-1}$ we get
$$
|t(x,y;u) - t(\xi,y;u)| \le cn^{-1}\rho(x, \xi)(1+n\rho(x, y)).
$$
This estimate coupled with \eqref{AxAy} leads to
\begin{equation}\label{Est-E1}
E_1 \le cn^d\rho(x, \xi)\frac{n^{2\mu}}{(1+n\rho(x, y))^\mu}
\le \frac{cn^d\rho(x, \xi)}{(A(x)+n^{-1})^\mu(A(y)+n^{-1})^\mu}.
\end{equation}

To estimate $E_2$ we observe that $t(x, y;u) = -t(x, -y;-u)$ and hence
$$
|t(x,y;u) - t(\xi,y;u)|= |t(x,-y;-u) - t(\xi,-y;-u)|.
$$
Consequently, $E_2$ can be estimated exactly as $E_1$ with
the same bound as in \eqref{Est-E1}. These two estimates yield \eqref{ball-Lip}.

\medskip


{\it Case 2:} $A(x)A(y) >  8c^*(n^{-1} \rho(x,y) + n^{-2})$.
In this case by \eqref{ball-1} and \eqref{Ax-Ay} it readily follows that
\begin{equation} \label{ball-AA}
A(x)A(y) \sim (A(x) +n^{-1})(A(x) +n^{-1}).
\end{equation}

Let $\eta_+$ be a $C^\infty$ function such that $\eta_+(u) = 1$ for $\frac12 \le u \le 1$,
and $\eta_+(u) = 0$ for $-1 \le u \le - \frac12$. Define $\eta_-(u) := 1-\eta_+(u)$.
Then on account of \eqref{ball-reprod-kernel}, we can write
\begin{align*}
  P_n(w_\mu; x,y) = P_n^+(w_\mu; x,y) + P_n^-(w_\mu; x,y),
\end{align*}
where
$$
  P_n^\pm (w_\mu;x,y) : =   c_\mu \frac{n+\lambda}{\lambda}
     \int_{-1}^1 C_n^\lambda(t(x,y;u)) \eta_\pm (u) (1-u^2)^{\mu-1} du.
$$
Since $t(x,y; - u) = -t(x, -y; u)$ and $C_n^\lambda(-t) = (-1)^n C_n^\lambda(t)$,
we only need to prove \eqref{ball-Lip} for $P_n^+(w_\mu;\cdot,\cdot)$.
We write $t(x,y;u)$ as
$$
  t(x,y;u) = B(x,y) + A(x)A(y)(u-1) \quad \hbox{with} \quad B(x,y) := \cos \rho (x,y),
$$
In going further, we have
\begin{equation}\label{ball-Lip+}
 P_n^+(w_\mu; x,y) - P_n^+(w_\mu; \xi,y) = J_1 + J_2,
\end{equation}
where
\begin{align*}
  J_1 := & c_\mu \frac{n+\lambda}{\lambda}
     \int_{-1}^1 \left[C_n^\lambda\big(B(x,y)+A(x)A(y) (u-1)\big) \right. \\
  &  \qquad\qquad\qquad  \left.  - C_n^\lambda\big(B(\xi,y)+A(x)A(y)(u-1)\big) \right]
            \eta_+ (u) (1-u^2)^{\mu-1} du, \\
  J_2 : = & c_\mu \frac{n+\lambda}{\lambda}
     \int_{-1}^1 \left[C_n^\lambda\big(B(\xi,y)+A(x)A(y) (u-1) \big) \right. \\
  &  \qquad\qquad\qquad  \left.  - C_n^\lambda\big(B(\xi,y)+A(\xi)A(y) (u-1)\big) \right]
            \eta_+(u) (1-u^2)^{\mu-1} du.
\end{align*}
To estimate  $|J_1|$, we again use
$\frac{d}{ds} C_n^\lambda(s) = 2 \lambda C_{n-1}^{\lambda+1} (s)$ to write
$$
J_1 = 2 c_\mu (n+\lambda) \int_{B(\xi,y)}^{B(x,y)}\int_{-1}^1
C_{n-1}^{\lambda+1}\big(s + A(x)A(y)(u-1)\big)\eta_+(u)(1-u^2)^{\mu-1} du  ds.
$$
We estimate the inner integral above using Lemma~\ref{ball-lem2} with $\eta(t)=\eta_+(t)(1+t)^{\mu-1}$,
$b = s-A(x)A(y)$ and $a = A(x) A(y)$. We get
\begin{align*}
 |J_1| \le c \frac{n^{2 (\lambda+1) - 2 \mu} }{[A(x)A(y)]^{\mu}}
 \Big|\int_{B(\xi,y)}^{B(x,y)} \frac{1}{(1+ n \sqrt{1 - |s|})^{\lambda+1-\mu}} ds \Big|.
\end{align*}
As in \eqref{t-t}
\begin{align*}
|B(x,y) - B(\xi,y)| \le \rho(x, \xi)(\rho(x, y)+\rho(\xi, y))
\le c \rho(x,\xi)(\rho (x,y) + n^{-1}).
\end{align*}
On the other hand
$
1-B(x, y) = 1-\cos\rho(x,y) \ge c\rho(x, y)^2
$
and similarly
$1-B(\xi, y) \ge c\rho(\xi, y)^2$.
Therefore,
\begin{equation}\label{est-J1}
|J_1| \le c  \frac{n^{ d+1} }{[A(x)A(y)]^{\mu}}  \frac{ (\rho (x,y) + n^{-1}) \rho(x,\xi)}
{\left(1+  n \min\{\rho(x, y), \rho(\xi, y)\} \right)^{\lambda +1 - \mu} }
\le \frac{cn^d \rho(x,\xi)}{[A(x)A(y)]^{\mu}},
\end{equation}
where we used that
$\rho(x, y)\le \rho(\xi, y)+\rho(x, \xi) \le \rho(\xi, y)+ c^*n^{-1}$.


To estimate $|J_2|$, we again use  $\frac{d}{ds} C_n^\lambda(s) = 2 \lambda C_{n-1}^{\lambda+1}(s)$
to express $J_2$ as
$$
-2 c_\mu (n+\lambda) A(y)\int_{A(\xi)}^{A(x)}\int_{-1}^1
C_{n-1}^{\lambda+1}\big(B(\xi,y) + sA(y) (u-1)\big)\eta_+(u)(1+u)^{\mu-1}(1-u)^{\mu} du  ds.
$$
We estimate the inner integral by using Lemma \ref{ball-lem2} with $\eta(t)=\eta_+(t)(1+t)^{\mu-1}$,
$b = B(\xi,y)-sA(y)$, $a = s A(y)$, and $\lambda$, $\mu$ replaced by $\lambda+1$, $\mu +1$ to obtain
\begin{align*}
 |J_2| \le \frac{cn^{2 (\lambda+1) - 2 (\mu+1)} }{A(y)^{\mu}{\big(1+ n \sqrt{1 - B(\xi,y)}\big)^{\lambda-\mu}}}
 \Big|\int_{A(\xi)}^{A(x)} s^{-\mu-1} ds \Big|.
\end{align*}
and using that $|A(x)-A(\xi)| \le \sqrt{2}\rho(x, \xi)$ (see \eqref{Ax-Ay})
\begin{equation}\label{est-J2}
|J_2| \le  \frac{cn^{d-1} \rho(x,\xi)}{A(y)^{\mu} \min\{A(x)^{\mu +1}, A(\xi)^{\mu+1}\}}.
\end{equation}
By the same token and since by assumption $\rho(x, \xi) \le c^*n^{-1}$ we have
$$
A(\xi) \ge A(x) -|A(x)-A(\xi)| \ge A(x) - \sqrt{2}\rho(x, \xi) \ge A(x)- \sqrt{2}c^*n^{-1}.
$$
If $A(x) \ge 2\sqrt{2}c^*n^{-1}$, then from above $A(\xi) \ge A(x)/2$.
These two estimates and \eqref{est-J2} imply that $|J_2|$ has the bound of $|J_1|$ from \eqref{est-J1},
and using \eqref{ball-AA} estimate \eqref{ball-Lip} holds for
$|P_n^+(w_\mu; x,y) - P_n^+(w_\mu; \xi,y)|$.

Let $A(x) < 2\sqrt{2}c^*n^{-1}$.
We claim that $A(y) < 4\sqrt{2}c^* n^{-1}$.
Indeed, suppose $A(y) \ge 4\sqrt{2}c^* n^{-1}$.
Then $A(y)\ge A(x)/2$ and using \eqref{Ax-Ay}, we get
$$
\sqrt{2}\rho(x, y) \ge |A(x)-A(y)| \ge A(y)-A(x) \ge A(y)/2
$$
and hence
$A(x)A(y) > 8c^*n^{-1}\rho(x, y)\ge 2\sqrt{2}c^*n^{-1}A(y)$
yielding $A(x) > 2\sqrt{2}c^*n^{-1}$,
that is a contradiction.
Therefore,
$A(x)A(y) < 16c^*n^{-2}$.
Thus $A(x)$, $A(y)$ obey the conditions of Case 1 and hence estimate \eqref{ball-Lip}
holds true.
This complete the proof of Lemma~\ref{ball-lem3}.
$\qed$

\end{document}